\documentclass[review,10pt]{elsarticle}

\usepackage{amsmath} 
\usepackage{amssymb}
\usepackage{amsfonts}
\usepackage{amsthm}
\usepackage{xcolor}
\usepackage{graphicx}
\usepackage{enumerate}
\usepackage{epstopdf} %converting to PDF
\usepackage{array,multirow}
\usepackage{tikz}
\usetikzlibrary{shapes, shadows, arrows,decorations.pathreplacing,calligraphy}
\tikzstyle{block}=[draw,rectangle,fill=blue!5,text width=12 em,text centered, minimum height=12mm, node distance=5 em]
\usepackage{caption}
\usepackage{subcaption}
\tikzstyle{line} = [draw,-latex']
\usetikzlibrary{calc}
\newcommand{\mysquare}{\rule[0.0mm]{2.2mm}{2.2mm}~}

\pgfkeys{%
	/polargrid/.cd,
	rmin/.code ={\global\def\rmin {#1}},
	rmax/.code ={\global\def\rmax {#1}},
	amin/.code ={\global\def\amin {#1}},
	amax/.code ={\global\def\amax {#1}},
	rstep/.code={\global\def\rstep{#1}}, 
	astep/.code={\global\def\astep{#1}}}
\def\polargrid{\pgfutil@ifnextchar[{\polar@grid}{\polar@grid[]}}%
\def\polargrid[#1]{%
	\pgfkeys{/polargrid/.cd,
		rmin ={2},
		rmax ={7},
		amin ={0},
		amax ={360},
		rstep={1}, 
		astep={22.5}}   
	\pgfqkeys{/polargrid}{#1}%
	\pgfmathsetmacro{\addastep}{\amin+\astep} 
	\pgfmathsetmacro{\addrstep}{\rmin+\rstep} 
	\foreach \a in {\amin,\addastep,...,\amax}  \draw[gray] (\a:\rmin) -- (\a:\rmax);  
	\foreach \r in {\rmin,\addrstep,...,\rmax}  \draw[gray] (\amin:\r cm) arc (\amin:\amax:\r cm);    
} 
\makeatother   
\usepackage{pgf}

\begin{document}

\begin{frontmatter}

%% Title, authors and addresses

%% use the tnoteref command within \title for footnotes;
%% use the tnotetext command for theassociated footnote;
%% use the fnref command within \author or \affiliation for footnotes;
%% use the fntext command for theassociated footnote;
%% use the corref command within \author for corresponding author footnotes;
%% use the cortext command for theassociated footnote;
%% use the ead command for the email address,
%% and the form \ead[url] for the home page:
%% \title{Title\tnoteref{label1}}
%% \tnotetext[label1]{}
%% \author{Name\corref{cor1}\fnref{label2}}
%% \ead{email address}
%% \ead[url]{home page}
%% \fntext[label2]{}
%% \cortext[cor1]{}
%% \affiliation{organization={},
%%             addressline={},
%%             city={},
%%             postcode={},
%%             state={},
%%             country={}}
%% \fntext[label3]{}

\title{A Generalised Curvilinear Coordinate system-based Patch Dynamics Scheme in Equation-free Multiscale Modelling}

%% use optional labels to link authors explicitly to addresses:
%% \author[label1,label2]{}
%% \affiliation[label1]{organization={},
%%             addressline={},
%%             city={},
%%             postcode={},
%%             state={},
%%             country={}}
%%
%% \affiliation[label2]{organization={},
%%             addressline={},
%%             city={},
%%             postcode={},
%%             state={},
%%             country={}}

\author[add1]{Tanay Kumar Karmakar \corref{cor1}} %% Author name
\address[add1]{Department of Mathematics, Indian Institute of Technology Guwahati, Assam-781039, India}
\ead{tanay.kumar@iitg.ac.in}
\cortext[cor1]{Corresponding author}

\author[add1]{Durga Charan Dalal}
\ead{durga@iitg.ac.in}

%% Author affiliation
%\affiliation{organization={Department of Mathematics, Indian Institute of Technology Guwahati},%Department and Organization
%            addressline={}, 
%            city={Guwahati},
%            postcode={781039}, 
%            state={Assam},
%            country={India}}

%% Abstract
\begin{abstract}
The patch dynamics scheme in equation-free multiscale modelling has the potential to efficiently predict the macroscopic behaviours by simulating the microscale problem in a fraction of the space-time domain. The patch dynamics schemes developed so far are mainly on rectangular domains with uniform grids and uniform rectangular patches. In real-life problems, the geometry of the domain is not regular or simple, where rectangular and uniform grids or patches may not be useful. To address this kind of complexity, for the first time, a generalised orthogonal curvilinear coordinate system is employed in the patch dynamics scheme, applicable to both rectangular domains with non-uniform grids and non-rectangular domains; while applying this, the concept of non-uniform and non-rectangular patch configurations in the physical domain is also adopted for the first time. An explicit representation of a patch dynamics scheme on a generalised curvilinear coordinate system in a two-dimensional domain is proposed for unsteady, linear, heterogeneous convection-diffusion-reaction ($\mathsf{CDR}$) problems. The robustness of the scheme on the generalised curvilinear coordinate system is assessed through numerical test cases. Firstly, a convection-dominated heterogeneous $\mathsf{CDR}$ equation is considered, featuring boundary layer regions in some part of the domain, for which stretched grids with non-uniform patch sizes are employed. The heterogeneity of the problems arises from the space-dependent diffusion tensor, convective velocity and the space-time-dependent source term. Spatially variable diffusivity is addressed in the first problem with spatially variable convective velocity. The results demonstrate that the non-uniform grid provides improved accuracy compared to the uniform grid in the presence of boundary layers. Secondly, a non-axisymmetric diffusion equation is examined in an annulus region, where the patches have non-rectangular shapes. The results obtained demonstrate excellent agreement with the analytical solution or existing numerical solutions. Overall, considering accuracy, computational time and memory usage, the proposed patch dynamics scheme performs much better compared to the full-domain simulation.
\end{abstract}

%%Graphical abstract
%\begin{graphicalabstract}
%\includegraphics{grabs}
%\end{graphicalabstract}

%%Research highlights
%\begin{highlights}
%\item A non-rectangular domain or a non-uniform grid can be handled using patch dynamics.
%\item Patch size and shapes can be non-uniform and non-rectangular.
%\item A patch dynamics scheme is proposed for generalized curvilinear coordinates.
%\item The technique is applied to unsteady linear convection-diffusion-reaction equations.
%\item Patch dynamics scheme can efficiently handle periodic boundary conditions.
%\end{highlights}

%% Keywords
\begin{keyword}
	multiscale modelling \sep equation-free framework \sep coarse projective integration \sep gap-tooth scheme \sep patch dynamics scheme
%multiscale modelling \sep equation-free framework \sep coarse integration \sep gap-tooth scheme \sep patch dynamics \sep curvilinear coordinate system

\end{keyword}

\end{frontmatter}

%\linenumbers

%%%%%%%%%%%%%%%%%%%%%%%%%%%%%%%%%%%%%%%%%%%%%%%%%%%%%%%%%%%%%%%%%%%%%%%%%%

\section{Introduction}
In order to capture the behaviours of a system accurately, the standard numerical schemes often require fine-level simulations. However, this fine-level simulation becomes futile for a large domain and a long duration due to the requirement of long computational time and huge memory. So, a balanced model that is computationally feasible and does not lose much information about the system is required; such a class of models may be acquired through multiscale modelling. %
Complex system \cite{2007_roberts_general,2014_roberts_dynamical} behaviours are generally characterised as coherent spatio-temporal dynamics arising from the interaction of its different components, and the classical numerical concepts mostly fail to handle such systems, whereas multiscale modelling becomes helpful. Almost all physical problems have multiple scales. Multiscale modelling captures the interactions and processes that occur at different scales in complex systems. In multiscale modelling, mostly the macroscopic and microscopic scales are used to understand the complex system behaviour. However, it becomes a challenge to couple the scales to exchange information between these scales.
In the fields of Science and Technology, macroscopic models always may not be available, though its microscopic models are available with sufficient accuracy \cite{1985_car_unified,1991_coron_numerical,1996_tadmor_quasicontinuum,2001_knap_analysis,2013_ortega_phenylalanine,2014_kiuchi_high,2014_nguyen_rupture}. Equation-free multiscale modelling \cite{2003_kevrekidis_equation} stands out as a powerful tool to take care of such problems. %Under the equation-free framework, one can handle $\mathsf{ODEs}$ or $\mathsf{PDEs}$ at the microscopic level to find the system-level behaviours. 

The time integration over a long duration in the equation-free framework proposed by Gear et al. \cite{2003_gear_projective} is known as projective integration ($\mathsf{PI}$). The forward Euler method and higher-order explicit method based on the Taylor series were used as outer integrators in $\mathsf{PI}$. An iterated projective method was also introduced by Gear et al. \cite{2003_gear_telescopic} for stiff problems, named as Telescopic Projective Methods. To extrapolate the macroscopic solution for a long time-step, forward Euler's method was used by some researchers \cite{2003_kevrekidis_equation,2005_samaey_gap,2006_samaey_patch,2020_arbabi_linking,2015_liu_acceleration,2009_samaey_equation,2003_gear_projective}. Lee et al. \cite{2007_lee_second} proposed second-order accurate Runge--Kutta and Adams--Bashforth methods as outer integrators to solve stiff differential systems.
Lafitte et al. \cite{2016_lafitte_high} and Maclean et al. \cite{2015_maclean_convergence} proposed higher-order projective integration schemes of arbitrary order based on the Runge--Kutta framework. These methods are designed specifically to address multiple time scales.

For a system having both spatial and temporal variables, the gap-tooth scheme in particle simulations was proposed by Gear et al. \cite{2003_gear_gap}. Kevrekidis et al. \cite{2003_kevrekidis_equation}  proposed the gap-tooth scheme ($\mathsf{GTS}$) to solve the evolution equation at the microscopic level for a short time in a large spatial domain. 
A more advanced algorithm in the equation-free multiscale modelling, which is a combination of gap-tooth scheme and coarse projective integration, is called the patch dynamics scheme, proposed by Kevrekidis et al. \cite{2003_kevrekidis_equation,2009_kevrekidis_equation} to find the system-level behaviours in a one-dimensional large domain and over a long time interval. A unique feature of the patch dynamics scheme is its capability to extrapolate macroscopic level behaviours from the microscopic level simulations on fractions of the space-time domain. This reduces the computational complexity and bridges intricate details and overarching patterns within complex systems \cite{2005_hyman_patch}. A practical introduction to the patch dynamics scheme is available in the Equation-Free Toolbox \cite{2021_maclean_toolbox}. Recently, Karmakar et al. \cite{2026_KUMARKARMAKAR_GPD} proposed generalised patch dynamics ($\mathsf{GPD}$) schemes in the equation-free framework, offering flexibility to utilise three distinct time scales: macro, meso and micro. The $\mathsf{GPD}$ scheme is of two types, $\mathsf{GPD-I}$ and $\mathsf{GPD-II}$. The $\mathsf{GPD}$ schemes are more robust and take less computational time compared to the other patch dynamics schemes. 

A $\mathsf{PDE}$ is called heterogeneous if one or more of its coefficients depend on its independent variables. Heterogeneous problems have been studied within the equation-free framework by Samaey et al. \cite{2005_samaey_gap,2006_samaey_patch}, Bunder et al. \cite{2017_bunder_good}, Karmakar et al. \cite{2026_KUMARKARMAKAR_GPD}.
%To achieve a better accuracy within less computational time compared to the usual patch dynamics schemes, one may use the $\mathsf{GPD}$ scheme.
%To decide the nature of the coarse system through microscopic simulation, like (i) the highest spatial derivative in the coarse equation, (ii) the dynamics of the coarse system is Hamiltonian or dissipative, or (iii) the coarse equation satisfies conservation laws or not, Li et al. \cite{2007_li_deciding} proposed ``The Baby-Bathwater Scheme". In these articles, the patch dynamics schemes were discussed only on one-dimensional domain.

So far, we have discussed patch dynamics schemes in one dimension. One of the initial works, as far as solving two-dimensional problems under the equation-free framework is concerned, Roberts et al. \cite{2014_roberts_dynamical} solved reaction-diffusion problems.  % and one of its applications is shown on the Ginzburg-Landau $\mathsf{PDE}$. 
%The concept of dynamical system theory is used to explore the macroscale modelling of PDEs as well as the microscale systems. The center manifold theory ensures the existence of a closed model on the macroscale grids. The ideas of the overlapping patches are also discussed. 
Bunder et al. \cite{2017_bunder_good} presented a two-dimensional patch dynamics scheme to solve problems involving discrete diffusion with fine-scale heterogeneity. 
%This article determines the optimum patch parameters from the patch dynamics simulation. 
They also showed the consistency of the macroscale dynamics with the closed microscale model in a self-adjoint patch dynamics scheme that provides an efficient, accurate and flexible computational homogenization \cite{2021_bunder_equation}. 

In two-dimensional space, gap-tooth schemes or patch dynamics schemes developed so far are mostly based on rectangular domains with uniform rectangular patches. However, the domain always may not be a rectangle for real-life problems. For the high-gradient regions in a domain, one needs to choose finer grids to capture the system-level behaviours accurately, that makes the grid spacings non-uniform. If patch dynamics schemes are used to solve such kind of problems, then the patch dynamics scheme based on non-uniform grids will be a better choice. Equation-free toolbox \cite{2020_Roberts_toolbox} allows non-uniform spacing of patches in one spatial dimension. Maclean et al. \cite{2021_maclean_equation,2022_maclean_adaptively} discussed the non-uniform size of patches, where they distributed the patches non-uniformly to adapt to resolving dynamically evolving shocks, particularly in a one-dimensional domain. For a non-rectangular macroscopic domain, body-fitted curvilinear coordinates can be implemented in such situations where the physical and geometrical aspects of a problem are better described using a new coordinate system. Some classical problems are solved using standard transformation techniques on non-uniform grids \cite{2007_pandit_transient,2008_pandit_fourth} and also on curvilinear coordinates \cite{2010_ray_transformation,2020_piquet_parallel}.

%In this article, the patch dynamics scheme \cite{2005_hyman_patch,2009_kevrekidis_equation,2020_arbabi_linking} is extended from one-dimensional space to two-dimensional space, and this patch dynamics scheme is used to solve general second-order unsteady convection-diffusion-reaction (CDR) equations at the microscopic level having nonuniform grids for the geometry of the domain that is beyond rectangular. Such models have inspired mathematical models for studying problems in many fields, such as fluid dynamics, heat transfer, physics of semiconductors, material engineering, chemistry, biology, population dynamics, astrophysics, biomedical engineering and financial mathematics. 

%Two-dimensional patch dynamics schemes developed till date have primarily focused on rectangular domains, where the patches are typically assumed to be uniform and rectangular in both shape and size. In contrast, 
This article proposes a new patch dynamics scheme to solve general second-order, unsteady, linear, heterogeneous convection-diffusion-reaction ($\mathsf{CDR}$) equations at the microscopic level on non-uniform grids and for complex domain geometries. Depending on the nature of the problem, the shapes and sizes of the patches may be chosen as either non-uniform or non-rectangular or both. To effectively handle such geometrically complex problems, generalised orthogonal curvilinear coordinates offer a natural and powerful framework. Accordingly, this article presents a two-dimensional patch dynamics scheme formulated on a generalised orthogonal curvilinear coordinate system in Section 3. To validate the proposed approach, two different types of numerical problems are considered based on their physical and geometrical complexities in Section 4. The first problem involves a two-dimensional advection-dominated heterogeneous $\mathsf{CDR}$ equation defined over a rectangular domain, characterised by steep gradients (or, boundary layer) near the right and top boundaries. This makes use of a stretched grid particularly advantageous in those regions. Both constant and variable diffusivity cases with variable convective velocity are addressed in the first problem.
%Non-homogeneous Dirichlet boundary conditions are considered, exhibiting exponential growth along the $x-$ and $y-$ directions. 
The second problem focuses on solving a two-dimensional, unsteady, non-axisymmetric diffusion equation in an annular domain using the body-fitted curvilinear coordinate system. %The numerical results obtained using the proposed scheme show excellent agreement with existing solutions, demonstrating the effectiveness and accuracy of the approach.

%%%%%%%%%%%%%%%%%%%%%%%%%%%%%%%%%%%%%%%%%%%%%%%%%%%%%%%%%%%%%%%%%%%%%%%%%%%%%%

\section{Governing equations in equation-free framework}
Let $u(x,y,t)$ denotes the microscopic state as a function of space variables $x$, $y$ and time variable $t$. We consider the unsteady microscopic problem in the two-dimensional domain $\Omega_p$ as 
\begin{equation}\label{eqn:Physical_Problem}
	\frac{\partial u}{\partial t}=\nabla\cdot\left(\mathcal{D}\nabla u-vu\right)+fu+g(x,y,t),
\end{equation}
where $\mathcal{D}=D(x,y)I_{2\times2}$ is the diffusion tensor, where $D(x,y)>0$ is the diffusion coefficient, $v(x,y)=\begin{bmatrix}
	v_1(x,y)\hspace{0.3cm} v_2(x,y)
\end{bmatrix}^\top$ is the non-negative convective velocity vector, $f$ is the reaction coefficient and $g(x,y,t)$ is the source term at point $(x,y)$ and at time $t$. Equation \eqref{eqn:Physical_Problem} represents a heterogeneous $\mathsf{CDR}$ equation, as both the diffusion coefficient and the convective velocity depend explicitly on the spatial variables and the source term depends on space-time variables.

\subsection{Evolution equation in a generalised orthogonal curvilinear form}

The development of a Cartesian coordinate system to a generalised orthogonal curvilinear coordinate system that facilitates computations based on the physical and geometrical aspects of the problem. This involves a transformation at the microscopic level from a complex physical domain to a rectangular computational domain with a uniform mesh. Thus, the simulation of the problem is performed on this computational domain. %Thus, the simulation formulating and solving the governing computational problem to perform the desired computations. 

Let the transformation be 
\begin{equation}\label{eqn:Transformtion}
	%	\begin{split}
		\xi=\xi(x,y), \hspace{0.2cm}\eta=\eta(x,y),
		%	\end{split}
\end{equation}
from the physical $xy$-plane to the computational $\xi\eta$-plane such that the transformation is non-singular. This transformation reshapes the intricated grids into elementary, uniform rectangular grids.

The corresponding inverse transformation of \eqref{eqn:Transformtion} may be written as
\begin{equation}\label{eqn:Transformtion_Inv}
	%	\begin{split}
		x=x(\xi,\eta), \hspace{0.2cm}y=y(\xi,\eta).
		%	\end{split}
\end{equation}

Under this transformation \eqref{eqn:Transformtion_Inv}, the microscopic problem \eqref{eqn:Physical_Problem} in the physical domain $\Omega_p$ becomes
\begin{equation}\label{eqn:Computational_Problem}
	\begin{split}
		\frac{\partial u}{\partial t}=\alpha(\xi,\eta)\frac{\partial^2u}{\partial \xi^2}+\beta(\xi,\eta)\frac{\partial^2 u}{\partial \xi \partial \eta}+\gamma(\xi,\eta)\frac{\partial^2u}{\partial \eta^2}\\
		-\nu(\xi,\eta)\frac{\partial u}{\partial \xi}-\omega(\xi,\eta)\frac{\partial u}{\partial \eta}
		+\phi(\xi,\eta) u+g(\xi,\eta,t),
	\end{split}
\end{equation}
in the computational domain $\Omega_c$, where
\begin{equation}
	\begin{split}
		\alpha(\xi,\eta)=\frac{D}{J^2}\left(x_\eta^2+y_\eta^2\right), \hspace{0.2cm} \beta(\xi,\eta)=-\frac{2D}{J^2}\left(x_\xi x_\eta+y_\xi y_\eta\right),\hspace{0.2cm}
		\gamma(\xi,\eta)=\frac{D}{J^2}\left(x_\xi^2+y_\xi^2\right), \hspace{0.2cm}\\
		\nu(\xi,\eta)=\frac{1}{J}\left[\left(x_\eta^2+y_\eta^2\right)\frac{\partial}{\partial\xi}\left(\frac{D}{J}\right)-\left(x_\xi x_\eta+y_\xi y_\eta\right)\frac{\partial}{\partial\eta}\left(\frac{D}{J}\right)\right]
		+\frac{D}{J^2}[\left(x_\eta x_{\eta\xi}-x_\xi x_{\eta\eta}\right)\\
		+\left(y_\eta y_{\eta\xi}-y_\xi y_{\eta\eta}\right)]-\frac{1}{J}\left[y_\eta v_1-x_\eta v_2\right],\\
		\omega(\xi,\eta)=\frac{1}{J}\left[\left(x_\xi^2+y_\xi^2\right)\frac{\partial}{\partial\eta}\left(\frac{D}{J}\right)-\left(x_\xi x_\eta+y_\xi y_\eta\right)\frac{\partial}{\partial\xi}\left(\frac{D}{J}\right)\right]
		+\frac{D}{J^2}[\left(x_\xi x_{\xi\eta}-x_\eta x_{\xi\xi}\right)\\
		+\left(y_\xi y_{\xi\eta}-y_\eta y_{\xi\xi}\right)]-\frac{1}{J}\left[x_\xi v_2-y_\xi v_1\right],\\
		\phi(\xi,\eta)=f-\frac{1}{J}\left[\left(y_\eta(v_1)_\xi-y_\xi(v_1)_\eta\right)+\left(x_\xi(v_2)_\eta-x_\eta(v_2)_\xi\right)\right], \hspace{0.2cm} J=x_\xi y_\eta-x_\eta y_\xi. 
	\end{split}
\end{equation}
Here, the notations denote $p_q=\frac{\partial p}{\partial q}$, $p_{qr}=\frac{\partial^2 p}{\partial r\partial q}$, where $p\in\{x,y,v_1,v_2,\frac{D}{J}\}$, $q,r\in\{\xi,\eta\}$. $\psi(\xi,\eta)$ means $\psi$ is evaluated at ($x(\xi,\eta)$, $y(\xi,\eta)$), where $\psi\in\{D,v_1,v_2\}$ and $g(\xi,\eta,t)$ means $g$ is evaluated at ($x(\xi,\eta)$, $y(\xi,\eta)$, $t$). The computational problem \eqref{eqn:Computational_Problem} corresponds to a heterogeneous $\mathsf{CDR}$ equation, as its coefficients depend on the spatial variables $\xi$ and $\eta$. In this article, we consider a generalised orthogonal curvilinear ($\xi-\eta$) coordinate system such that the coefficient of the mixed derivative (arising due to the transformation) vanishes, i.e., $x_\xi x_\eta+y_\xi y_\eta$=0. 

The curvilinear transformation \eqref{eqn:Transformtion_Inv} is constructed in such a way that it is smooth and one-to-one over the entire computational domain. The Jacobian determinant $J=x_\xi y_\eta-x_\eta y_\xi$ remains strictly positive and bounded throughout the domain. These guarantee local invertibility and prevent singularities in the transformed $\mathsf{PDE}$ \cite{1998_thompson_handbook,1999_liseikin_grid,2003_farrashkhalvat_basic}. In addition, the mapping we use is free from grid folding or overlapping, ensuring global invertibility and well-behaved metric coefficients. The transformation \eqref{eqn:Transformtion_Inv} does not introduce singularities, and the $\mathsf{PDE}$ \eqref{eqn:Computational_Problem} remains well-posed in the ($\xi,\eta$) coordinate system.

\section{Patch dynamics scheme in two-dimensional space}

Using the transformation \eqref{eqn:Transformtion_Inv}, the microscopic problem \eqref{eqn:Physical_Problem} is transformed into a computational problem \eqref{eqn:Computational_Problem} on the rectangular domain $\Omega_c=[a,b]\times [c,d]$.
A patch dynamics scheme is implemented on the domain $\Omega_c$ to find system-level behaviours of the computational problem \eqref{eqn:Computational_Problem}.
We discretise the macroscopic space domain $\Omega_c=[a,b]\times[c,d]$ of the computational problem using equidistant macroscopic mesh $\{(\xi_i, \eta_j) \mid 0 \le i \le N_\xi,\ 0 \le j \le N_\eta,\ \xi_i = a + i\Delta \xi,\ \eta_j = c + j\Delta \eta\}$, ensuring that a closed approximation is achieved, where $\Delta\xi = \frac{b - a}{N_\xi}, \quad \Delta\eta = \frac{d - c}{N_\eta}$ are the macroscopic grid sizes along $\xi$- and $\eta$-directions, respectively. In the computational domain $\Omega_c$, the patches are rectangular in shape and uniform in size. We define the length of the edges of the patches are as $h_\xi$ and $h_\eta$ along $\xi$- and $\eta$-axes, respectively, and centered about the macroscopic grid ($\xi_{i},\eta_{j}$) for all $i=1,\ldots,(\operatorname{N_\xi-1}),\hspace{0.2cm} j=1,\ldots,(\operatorname{N_\eta-1})$. We choose $h_\xi$ and $h_\eta$ in such a way that they satisfy $\frac{h_\xi}{\Delta\xi}<1$ and $\frac{h_\eta}{\Delta\eta}<1$, this also ensures that the patches do not overlap. The microscopic problem \eqref{eqn:Computational_Problem} is computed within each of the patches in the computational domain. For that, we discretise the microscopic problem using a finite difference scheme on each patch for a micro time step $\tau$, where $\tau$ satisfies the relaxation time limit. Suppose, $\delta \xi$, $\delta \eta$ are the spatial nano-steps along $\xi$- and $\eta$-directions, respectively, and $\delta t$ is the temporal nano-step. Using averaging of the microscale solution over the $(i,j)^{th}$ patch, we restrict a macroscale solution $U_{i,j}(t)$ at ($\xi_{i}$,$\eta_{j}$) by
\begin{equation}\label{eqn:Averaging}
	{U}_{i,j}(t):=\frac{1}{h_\xi h_\eta}\int_{\eta_{j}-\frac{h_\eta}{2}}^{\eta_{j}+\frac{h_\eta}{2}}\int_{\xi_{i}-\frac{h_\xi}{2}}^{\xi_{i}+\frac{h_\xi}{2}}u(w,z,t) \hspace{0.1cm}dw \hspace{0.1cm}dz.
\end{equation} 
This restriction operator \eqref{eqn:Averaging} is a subjective preference for the user. The same restriction operator was also used by Bunder et al. \cite{2017_bunder_good,2021_bunder_equation}.

\subsection{Patch edge conditions}

It is essential to select appropriate patch edge conditions based on the nature of the problem in small patches $\Omega_c^{i,j}:=\left[\xi_i-\frac{h_\xi}{2},\xi_i+\frac{h_\xi}{2}\right]\times\left[\eta_j-\frac{h_\eta}{2},\eta_j+\frac{h_\eta}{2}\right]$, $\forall i=1,...,N_\xi-1$, $\forall j=1,...,N_\eta-1$. We approximate the macroscopic field $U$ using a bi-quadratic polynomial within $\Omega_c^{i,j}$ patch as:
\begin{equation}\label{eqn:Polynomial}
	u(\xi,\eta,t)\approx \mathcal{P}_{i,j}(\xi,\eta,t), \hspace{0.15cm} (\xi,\eta)\in\Omega_c^{i,j},
\end{equation}
where $\mathcal{P}_{i,j}$ represents a Lagrange polynomial,
\begin{equation}\label{eqn:Lagrange_2D}
	\mathcal{P}_{i,j}(\xi,\eta,t)=\sum_{p=i-1}^{i+1}\sum_{q=j-1}^{j+1} L_p^i(\xi)L_q^j(\eta)U_{p,q}(t),
\end{equation}
within $(i,j)^{th}$ patch, where, 
\begin{equation}
	L_\beta^\gamma(\alpha)=\prod_{\begin{array}{c}
			r=\gamma-1 \\ r\neq \beta
	\end{array}}^{\gamma+1}\frac{\alpha-\alpha_r}{\alpha_\beta-\alpha_r}
\end{equation} 
are the Lagrange fundamental polynomials of degree 2, $\forall i=1,...,(N_\xi-1)$ and $\forall j=1,...,(N_\eta-1)$.  %$U_{p,q}(t)$ be the approximate value of the coarse variable $U$ at the macro grid ($\xi_p,\eta_q$) at time $t$.

\begin{figure}%[h!]
	\centering
	\begin{tikzpicture}[scale=1.5]
		\draw[step=1cm,gray, very thick] (-0.9,-0.9) grid (2.9,2.9);
		
		\draw[step=0.1cm, thin] (-0.2,-0.2) grid (0.2,0.2);
		\draw[thick] (-0.2,-0.2) rectangle (0.2,0.2);
		\draw[step=0.1cm, thin] (0.8,-0.2) grid (1.2,0.2);
		\draw[thick] (0.8,-0.2) rectangle (1.2,0.2);
		\draw[step=0.1cm, thin] (1.8,-0.2) grid (2.2,0.2);
		\draw[thick] (1.8,-0.2) rectangle (2.2,0.2);
		
		\draw[step=0.1cm, thin] (-0.2,0.8) grid (0.2,1.2);
		\draw[thick] (-0.2,0.8) rectangle (0.2,1.2);
		\draw[step=0.1cm, thin] (0.8,0.8) grid (1.2,1.2);
		\draw[thick] (0.8,0.8) rectangle (1.2,1.2);
		\draw[step=0.1cm, thin] (1.8,0.8) grid (2.2,1.2);
		\draw[thick] (1.8,0.8) rectangle (2.2,1.2);
		
		\draw[step=0.1cm, thin] (-0.2,1.8) grid (0.2,2.2);
		\draw[thick] (-0.2,1.8) rectangle (0.2,2.2);
		\draw[step=0.1cm, thin] (0.8,1.8) grid (1.2,2.2);
		\draw[thick] (0.8,1.8) rectangle (1.2,2.2);
		\draw[step=0.1cm, thin] (1.8,1.8) grid (2.2,2.2);
		\draw[thick] (1.8,1.8) rectangle (2.2,2.2);
		
		\filldraw[color=blue!100, fill=blue!100, very thick](0,0) circle (0.05);
		\filldraw[color=blue!100, fill=blue!100, very thick](1,0) circle (0.05);
		\filldraw[color=blue!100, fill=blue!100, very thick](2,0) circle (0.05);
		\filldraw[color=blue!100, fill=blue!100, very thick](0,1) circle (0.05);
		\filldraw[color=red!100, fill=red!100, very thick](1,1) circle (0.05);
		\filldraw[color=blue!100, fill=blue!100, very thick](2,1) circle (0.05);
		\filldraw[color=blue!100, fill=blue!100, very thick](0,2) circle (0.05);
		\filldraw[color=blue!100, fill=blue!100, very thick](1,2) circle (0.05);
		\filldraw[color=blue!100, fill=blue!100, very thick](2,2) circle (0.05);
		
		\draw[very thick] (1.35,-0.4) node {\small{$i,j-1$}};
		\draw[very thick] (0.2,-0.4) node {\small{$i-1,j-1$}};
		\draw[very thick] (2.4,-0.4) node {\small{$i+1,j-1$}};
		
		\draw[very thick] (1.3,0.6) node {\small{$i,j$}};
		\draw[very thick] (0.23,0.6) node {\small{$i-1,j$}};
		\draw[very thick] (2.35,0.6) node {\small{$i+1,j$}};
		
		\draw[very thick] (1.35,1.6) node {\small{$i,j+1$}};
		\draw[very thick] (0.075,1.6) node {\small{$i-1,j+1$}};
		\draw[very thick] (2.4,1.6) node {\small{$i+1,j+1$}};
		
	\end{tikzpicture}
	\vspace{0.2cm}
	\caption{A schematic view of the 9-point stencil used in equations \eqref{eqn:Polynomial}, \eqref{eqn:Lagrange_2D} and \eqref{eqn:edge}. The rectangles represent the patches where microscopic simulations are performed. To determine the patch edge conditions for the $(i,j)^{th}$ patch, we communicate with neighbouring macroscopic grid points (indicated by blue dots), including the macro grid point at the centre of the patch itself (indicated by a red dot).}
	\label{fig:Stencil_9_Point}
\end{figure}

Polynomial interpolation serves as a bridge across spatial gaps between patches, implementing communication between them. Equation \eqref{eqn:Lagrange_2D} serves as an interpolating polynomial between the neighbouring macroscopic grid points \{$U_{p,q}(t)$, $p\in \{i-1, i, i+1\}, q\in \{j-1, j, j+1\}$\} for the $(i,j)^{th}$ patch at a fixed time $t$. %We define a non-zero flux value for edge conditions to maintain continuous diffusion. In order to accomplish this, we use Neumann's condition. 
We propose Neumann edge conditions as the gradients of $u$ at the four edges of each patch as,
\begin{equation}\label{eqn:edge}
	\begin{split}
		\frac{\partial u}{\partial\xi}\bigg|_{\xi_i\pm \frac{h_\xi}{2}} = \frac{1}{\Delta\xi} \Bigg[ 
		& \left(-\frac{1}{2} \pm r_\xi \right) \sum_{q=j-1}^{j+1} L^j_q(\eta) U_{i-1,q}
		\mp 2r_\xi \sum_{q=j-1}^{j+1} L^j_q(\eta) U_{i,q} \\
		& + \left( \frac{1}{2} \pm r_\xi \right) \sum_{q=j-1}^{j+1} L^j_q(\eta) U_{i+1,q} 
		\Bigg],\\
		\frac{\partial u}{\partial\eta}\bigg|_{\eta_j\pm \frac{h_\eta}{2}} = \frac{1}{\Delta\eta} \Bigg[ 
		& \left(-\frac{1}{2} \pm r_\eta \right) \sum_{p=i-1}^{i+1} L^i_p(\xi) U_{p,j-1}
		\mp 2r_\eta \sum_{p=i-1}^{i+1} L^i_p(\xi) U_{p,j} \\
		& + \left( \frac{1}{2} \pm r_\eta \right) \sum_{p=i-1}^{i+1} L^i_p(\xi) U_{p,j+1} 
		\Bigg],
	\end{split}
\end{equation}
at time $t$. The length scale ratioes are denoted as $r_\xi:=\frac{h_\xi}{2\Delta\xi}$ and $r_\eta:=\frac{h_\eta}{2\Delta\eta}$.

We keep these derivative values fixed during micro-simulation for the entire duration of the time-stepper [$t,t+\tau$]. Keeping these values fixed is the $Q$=1 case of Bunder et al. \cite{2016_bunder_accuracy}. $Q>1$ is a more accurate alternative for future work. In the equation-free toolbox, Roberts et al. \cite{2020_Roberts_toolbox} demonstrated that the patch edge values are automatically and dynamically updated at each sub-step of the microscale integration, based on the user's chosen configuration. In large-scale problems, one may wish to limit communication between patches only at mesoscopic time scales, which are notably larger than the microscopic time scale. An approach to achieve this was proposed by Bunder et al. \cite{2016_bunder_accuracy}.

In order to construct the patch edge conditions for the $(i,j)^{th}$ patch in our proposed scheme, a total of 9 neighbouring macroscopic values, including the central value $U_{i,j}$ itself, are used to establish communication, as depicted in Figure \ref{fig:Stencil_9_Point}.

\subsection{Initial condition}
Suppose, the values of the macroscopic variable $U$ are either provided (at $t=0$) or computed at the macro grid points ($\xi_{i},\eta_{j}$), $\forall$$i, j$ at time $t$.
For time integration, initial values of $u_{i,j}(\xi,\eta,t)$ are required on each patch $\Omega_c^{i,j}$ at time $t$. The initial values of $u_{i,j}(\xi,\eta,t)$ are obtained from $U(\xi,\eta,t)$ through a lifting operator in the equation-free framework as shown in equation \eqref{eqn:Lifting}. %As a standard choice, we assume the lifting operator as a polynomial expansion within the $(i,j)^{th}$ patch, expressed as:
\begin{equation}\label{eqn:Lifting}
	u_{i,j}(\xi,\eta,t)=C_0(t)+\sum_{l=1}^{2}\frac{1}{l!}\left((\xi-\xi_i)\frac{\partial}{\partial\xi}+(\eta-\eta_j)\frac{\partial}{\partial\eta}\right)^l U(\xi,\eta,t),
\end{equation}
where $(\xi,\eta)\in \Omega_c^{i,j}$ at time $t$. The term, $C_0$, is obtained from the micro variable $u$ using the patch averaging property \eqref{eqn:Averaging} as 
\begin{equation}\label{eqn:C_0}
	C_0(t)=\left[U-\frac{1}{24}h_\xi^2\frac{\partial^2U}{\partial \xi^2}-\frac{1}{24}h_\eta^2\frac{\partial^2U}{\partial \eta^2}\right]\bigg|_{(\xi_{i},\eta_{j},t)}.
\end{equation}
All partial derivative terms in \eqref{eqn:Lifting} and \eqref{eqn:C_0} are discretised using proper finite difference schemes at the grid location $(\xi_i,\eta_j)$ at time $t$. 

\subsection{Algorithm of patch dynamics ($\mathsf{PD}$) scheme}
Patch dynamics scheme replaces expensive nano-scale simulations across a large domain over a long time interval with computations across small, sparse patches in space-time. The patch dynamics ($\mathsf{PD}$) scheme is a combination of the gap-tooth scheme and coarse projective integration. Consider that the entire time interval $[0, T]$ is discretised with a macroscopic time step $\Delta t=\frac{T}{\operatorname{N_t}}$, where $\operatorname{N_t}$ be the total macro time levels. Here, $n\Delta t$ is defined as $T_n$ for $0\le n\le\operatorname{N_t}$. Suppose at time $T_n$, the values of $U_{i,j}(T_n)$ at the macro grid points $(\xi_{i},\eta_{j})$, $\forall$ $i,j$ are available. Below, we provide a complete algorithm to progress a macro time step $\Delta t$ from one macro time level $T_n$ to the next macro time level $T_{n+1}$ in the computational domain:

\begin{itemize}
	\item[\mysquare]\textbf{Gap-tooth scheme ($\mathsf{GTS}$)}
\end{itemize}
Consider that we know the values of $U_{i,j}^{0,n}$ at time $T_n$ for all macro grid points $(\xi_{i},\eta_{j})$, $\forall$ $i,j$. The notation $U_{i,j}^{p,q}$ represents the macroscale solution at the macroscopic grid points $(\xi_{i},\eta_{j})$ at time $q\Delta t+p\tau$, where $\Delta t$ is the macro time step, and $\tau$ is the micro evolution time for $\mathsf{GTS}$. We construct the $\mathsf{GTS}$ using the following steps:
\begin{enumerate}[(i)]	
	
	\item \textbf{Edge conditions:} We compute the values required for the patch edge conditions using equation \eqref{eqn:edge} at time $T_n$. We keep these values fixed during micro-simulation for the short time interval [$T_n, T_n+\tau$]. 
	\item \textbf{Lifting:} At time $T_n$, we formulate the microscopic initial condition $u_{i,j}(\xi,\eta,T_n)$ inside the patch using \eqref{eqn:Lifting}. 
	
	\item \textbf{Evolution:} 
	We compute the microscopic problem \eqref{eqn:Computational_Problem} inside the patches $\Omega_c^{i,j}$, incorporating the initial condition \eqref{eqn:Lifting}, and the edge conditions \eqref{eqn:edge} until time $T_n+\tau$. 
	
	\item \textbf{Restriction:}
	We compute a macroscale value
	\begin{equation}\label{eqn:Restriction}
		U_{i,j}^{1,n}=\frac{1}{h_\xi h_\eta}\int_{\eta_{j}-\frac{h_\eta}{2}}^{\eta_{j}+\frac{h_\eta}{2}}\int_{\xi_{i}-\frac{h_\xi}{2}}^{\xi_{i}+\frac{h_\xi}{2}}u_{i,j}(w,z,t_n+\tau)dw\hspace{0.1cm}dz,
	\end{equation}	
	at macro grid point $(\xi_i,\eta_{j})$ at time $T_n+\tau$ using averaging technique over the micro solution $u$ in the $(i,j)^{th}$ patch at time $T_n+\tau$ for all $i=1,...,(N_\xi-1) \hspace{0.25cm} \text{and}\hspace{0.25cm} j=1,...,(N_\eta-1)$. %, via a quadrature formula.

	\begin{itemize}
		\item[\mysquare]\textbf{Coarse projective integration}
	\end{itemize}
	\item \textbf{ Short time step:} We perform gap-tooth scheme once to determine $U_{i,j}^{1,n}$ at time $T_n+\tau$ (points (i) through (iv)), aiding in the evaluation of an approximate value of the time derivative for the macroscopic variable, i.e., $\frac{\partial U}{\partial t}$ at $(\xi_i,\eta_j, T_n)$.
	
	The time derivative of the macro state in each patch at time $T_n$ is estimated using forward difference as:
	\begin{equation}\label{eqn:Time_Derivative}
		F(\xi_{i}, \eta_{j}, T_n):=\frac{U_{i,j}^{1,n}-U_{i,j}^{0, n}}{\tau}, \hspace{0.1cm} \forall i=1,..., (N_\xi-1) \hspace{0.15cm} \text{and}\hspace{0.15cm} j=1,..., (N_\eta-1).
	\end{equation}
	\item \textbf{Projective extrapolation:}
	We use the above estimate with the forward Euler method to advance $U$ forward over a time step $\Delta t$ to reach the next coarse level. For forward Euler, the extrapolation is given by: 
	\begin{equation}\label{eqn:Extrapolation}
		U_{i,j}^{0,n+1}=U_{i,j}^{0,n}+\Delta t\hspace{0.05cm} F(\xi_{i}, \eta_{j}, T_n), \hspace{0.1cm} \forall i=1,..., (N_\xi-1) \hspace{0.15cm} \text{and}\hspace{0.15cm} \forall j=1,..., (N_\eta-1).
	\end{equation}
	We repeat the entire procedure (points (i) through (vi)) within each macro time step $\Delta t$ until the final time $T$ is reached. This way, one obtains the coarse solution in a large domain over a long time with the help of the microscopic problem, where the computation is performed on a small portion of the space-time domain.
\end{enumerate}

As an initial version, the proposed scheme is developed to address complex systems, specifically heterogeneous problems with smooth, non-oscillatory coefficients. However, the scheme presented in this study may not be fully suitable for homogenisation problems with oscillatory coefficients due to the use of the restriction operator \eqref{eqn:Restriction} and the artificial edge conditions \eqref{eqn:edge}. In such settings, local averaging of the micro solutions over the entire patch introduces a resonance error, which originates from these artificial edge conditions, as documented in \cite{2011_gloria_reduction,2016_arjmand_time,2019_abdulle_exponential,2021_abdulle_parabolic,2023_abdulle_elliptic}. In order to mitigate this issue, it is necessary to employ either buffer regions \cite{2006_samaey_patch,2020_arbabi_linking} or action--core regions \cite{2016_bunder_accuracy,2017_bunder_good} or the concept of averaging kernel \cite{2018_arjmand_equation}, which effectively shield the microscale simulations from spurious oscillations near the patch edges due to artificial boundary conditions. The proposed scheme will be further developed to handle heterogeneous problems with highly oscillatory coefficients in future work.

We observe that for numerical test problems that do not exhibit oscillatory behaviours at the microscale, the standard averaging-based restriction operator described in this study works well. This approach was successfully used in several earlier equation-free studies, such as Kevrekidis et al. \cite{2003_kevrekidis_equation,2004_theodoropoulos_coarse,2009_kevrekidis_equation}, Samaey et al. \cite{2005_samaey_gap,2009_samaey_equation}, Xiu et al. \cite{2005_xiu_equation}, Li et al. \cite{2007_li_deciding}, Roberts et al. \cite{2010_roberts_choose}, Karmakar et al. \cite{2026_KUMARKARMAKAR_GPD}, etc. Therefore, for the heterogeneous problems with smooth, non-oscillatory coefficients considered in our manuscript, the results in the following section demonstrate that the proposed restriction operator is appropriate and reliable.

%%%%%%%%%%%%%%%%%%%%%%%%%%%%%%%%%%%%%%%%%%%%%%%%%%%%%%%%%%%%%%%%%%%%%%%%%%

\section{Results and discussion}
The proposed patch dynamics scheme is validated through two test cases. The first case involves a heterogeneous convection-diffusion-reaction ($\mathsf{CDR}$) equation over a two-dimensional stretched grid in Subsection 4.1, where convective velocity depends on space variables. In the first scenario, the diffusion coefficient is constant, while in the second scenario, the diffusion coefficient depends on space variables. Additionally, in both scenarios, the problems exhibit convection dominance. The second test case explores non-axisymmetric diffusion in a two-dimensional annulus in Subsection 4.2. 

\subsection{Heterogeneous convection-diffusion-reaction ($\mathsf{CDR}$) equation}

The convection-diffusion-reaction equation plays a crucial role in capturing the behaviours of the systems that exhibit multiple spatial and temporal scales. $\mathsf{CDR}$ equations include chemical and biochemical processes, environmental transport, drug delivery, tumour growth and angiogenesis, fluid dynamics and heat transfer, materials science, biological systems, combustion modelling, atmospheric and oceanic modelling, nanotechnology, etc.

\subsubsection{Patch dynamics ($\mathsf{PD}$) solution with constant diffusivity:}

%The two-dimensional unsteady convection-diffusion-reaction equation is described as:

$\textbf{Physical Problem:}$

In equation \eqref{eqn:Physical_Problem}, we set 
\begin{equation}\label{eqn:P1_P}
	D(x,y)=1,\hspace{0.2cm} v=\begin{bmatrix}
		10x \hspace{0.3cm}10y
	\end{bmatrix}^\top,\hspace{0.2cm} f=20,\hspace{0.2cm} g(x,y,t)=(10x+10y-1)e^{x+y+t},
\end{equation}
in the physical domain $\Omega_p=[0,1]\times [0,1]$ with final time $T=1$. The initial condition is $u(x,y,0)=e^{x+y}$ in $\Omega_p$. The boundary conditions of the four physical boundaries are $u(0,y,t)=e^{y+t}$, $u(1,y,t)=e^{1+y+t}$, $u(x,0,t)=e^{x+t}$, and $u(x,1,t)=e^{x+1+t}$ on $\partial\Omega_p$ and $t\in[0, 1]$, where $\Omega_p$ is the physical domain, and $\partial \Omega_p$ is the boundary of $\Omega_p$. The problem is heterogeneous, with heterogeneity arising from the space-dependent advection velocity and the space-time-dependent source term.

%Here, $u$ is the temperature, and the thermal diffusivity along the $x$ and $y$ directions is $D_x=1$ and $D_y=1$, respectively. $v_x=10x$ and $v_y=10y$ are the components of convective velocity in the domain $\Omega_p$ along the $x$ and $y$ directions, respectively. The physical domain is $\overline{\Omega}_p=[0,L]\times [0,H]$, where $L=H=1$. The source term is $S(x,y,t)=(10x+10y-1)e^{x+y+t}$.

The exact solution of the physical problem \eqref{eqn:Physical_Problem} is given by:
\begin{equation}\label{eqn:CDR_Phy_Exact}
	u_e(x,y,t)=e^{x+y+t} \hspace{0.25cm} \text{in} \hspace{0.25cm} \Omega_p\times [0,1].
\end{equation}

Spatial discretisation would transform the $\mathsf{PDE}$ \eqref{eqn:Physical_Problem} into a large system of stiff $\mathsf{ODEs}$  \cite{2007_lee_second}, characterised by widely separated time scales. The local reaction term induces rapid changes in temperature over time, in contrast to the comparatively slower diffusion process, a hallmark of multiscale behaviour. %The fast variation in local temperature (or concentration) due to the reaction term clearly reflects the existence of multiple time scales. 
Additionally, the exponential nature of the solution leads to steep gradients (boundary layers) near the right and top boundaries in the physical domain $\Omega_p$, introducing a spatial multiscale nature. To efficiently address these challenges, a multiscale technique on both spatial and temporal scales is required. In the following section, we apply our proposed scheme to deal with such problems, which involve multiple scales in both space and time.
Since 99\% of the physical domain is convection-dominated, we classify the problem \eqref{eqn:Physical_Problem} in Sub-section 4.1 as a convection-dominated problem. In Subsubsection 4.1.1, the problem \eqref{eqn:Physical_Problem} is solved using the patch dynamics scheme in the primitive variables $x$ and $y$, whereas in Subsubsection 4.1.2, it is solved in the transformed variables $\xi$ and $\eta$.

\textbf{$\bullet$ Patch dynamics ($\mathsf{PD}$) solution without transformation:}

\begin{table}%[h!]
%	\begin{center}
	\caption{We compare the percentage errors of the half boundary method ($\mathsf{HBM}$) solution and the patch dynamics ($\mathsf{PD}$) solution with respect to the exact solution \eqref{eqn:CDR_Phy_Exact} for the problem \eqref{eqn:Physical_Problem} at final time $T=1$. The comparison is carried out using various uniform macroscopic spatial and temporal discretisations. In the last four columns, both percentage errors are reported at four distinct spatial locations: (0.2,0.2), (0.4,0.4), (0.6,0.6), and (0.8,0.8).}
	\centering
		\scalebox{1}{
			\begin{tabular}{ |p{2cm}|p{3.5cm}|p{1.3cm}|p{1.3cm}|p{1.3cm}|p{1.3cm}|}
				\hline
				\multicolumn{2}{|c|}{} & \multicolumn{4}{c|}{Percentage errors at nodes at final time} \\
				\hline
				Grids and time levels& Method & \multicolumn{1}{c|}{$\left(0.2, 0.2\right)$} & \multicolumn{1}{c|}{$\left(0.4, 0.4\right)$} & \multicolumn{1}{c|}{$\left(0.6, 0.6\right)$} & \multicolumn{1}{c|}{$\left(0.8, 0.8\right)$} \\
				\hline
				$10\times10$, 100&$\mathsf{HBM}$ (Zhao et al. \cite{2022_zhao_half})&2.67&4.76&5.42&4.41\\
				$10\times10$, 1000&$\mathsf{PD}$ scheme (present) &1.93$\mathrm{e}{-2}$&5.11$\mathrm{e}{-2}$&7.40$\mathrm{e}{-2}$&7.00$\mathrm{e}{-2}$\\
				\hline
				$15\times15$, 100&$\mathsf{HBM}$ (Zhao et al. \cite{2022_zhao_half})& 1.86&3.33&3.76&2.94\\
				$15\times15$, 2200&$\mathsf{PD}$ scheme (present)&5.38$\mathrm{e}{-3}$&1.68$\mathrm{e}{-2}$& 2.59$\mathrm{e}{-2}$&2.52$\mathrm{e}{-2}$\\
				\hline
				$20\times20$, 100&$\mathsf{HBM}$ (Zhao et al. \cite{2022_zhao_half})& 1.44&2.58&2.91&2.25\\
				$20\times20$, 4000&$\mathsf{PD}$ scheme (present)&4.80$\mathrm{e}{-4}$&4.73$\mathrm{e}{-3}$&9.18$\mathrm{e}{-3}$
				&9.91$\mathrm{e}{-3}$\\
				\hline
				$25\times25$, 100&$\mathsf{HBM}$ (Zhao et al. \cite{2022_zhao_half})&1.18&2.12&2.39&1.84\\
				$25\times25$, 6200&$\mathsf{PD}$ scheme (present)&1.78$\mathrm{e}{-3}$& 8.29$\mathrm{e}{-4}$&1.46$\mathrm{e}{-3}$&2.92$\mathrm{e}{-3}$\\
				\hline
		\end{tabular}}
%	\end{center}
	\label{table:HBMvsPD}
\end{table}

In this article, the proposed patch dynamics scheme is used to solve the heterogeneous problem \eqref{eqn:Physical_Problem}, which was also numerically solved using the half boundary method ($\mathsf{HBM}$) of Zhao et al. \cite{2022_zhao_half}. 
%In Table \ref{table:HBMvsPD}, a comparison is made between the $\mathsf{HBM}$ solution, patch dynamics solution and the exact solution of the problem \eqref{eqn:CDR_Phy_1} at (0.2, 0.2), (0.4, 0.4), (0.6, 0.6) and (0.8, 0.8) locations for different sizes of grids. 
A uniform discretisation of the physical domain is considered to compare with the $\mathsf{HBM}$ results of \cite{2022_zhao_half}. In the patch dynamics scheme, the 2-point upwind scheme of first order and the 3-point central difference scheme of second order are used for the convective and diffusive terms, respectively, in micro-simulation. The alternating direction implicit (ADI) scheme is used for the microsimulation inside the patches. The edge-lengths of each patch are considered as $h_x=h_y=0.001$, and the time-stepper size as $\tau=1e-6$ within the domain $\Omega_p$. At the microscopic level, each patch is discretised using $\operatorname{n}=10$ spatial nano grid steps along each $x$- and $y$-direction. The time within the time-stepper is discretised using $\operatorname{n_t}=2$ nano time steps. %Micro means patch scale, and nano means sub-patch discretisation.
Although the ADI scheme is an implicit scheme at the microscopic level, the patch dynamics scheme is explicit at the macroscopic level. To fulfil the stability criteria of the patch dynamics scheme at the macroscopic level, we need to take more number of macro-time steps than that of the $\mathsf{HBM}$ scheme. However, both schemes have the same macroscopic spatial grids. The trapezoidal composite rule is used in the restriction operator \eqref{eqn:Restriction} to restrict the microscopic values of the solution in the patch to its macroscopic value.

\begin{figure}%[h!]
	\begin{subfigure}{.45\textwidth}
		\centering
		% include first image
		\includegraphics[width=1\linewidth]{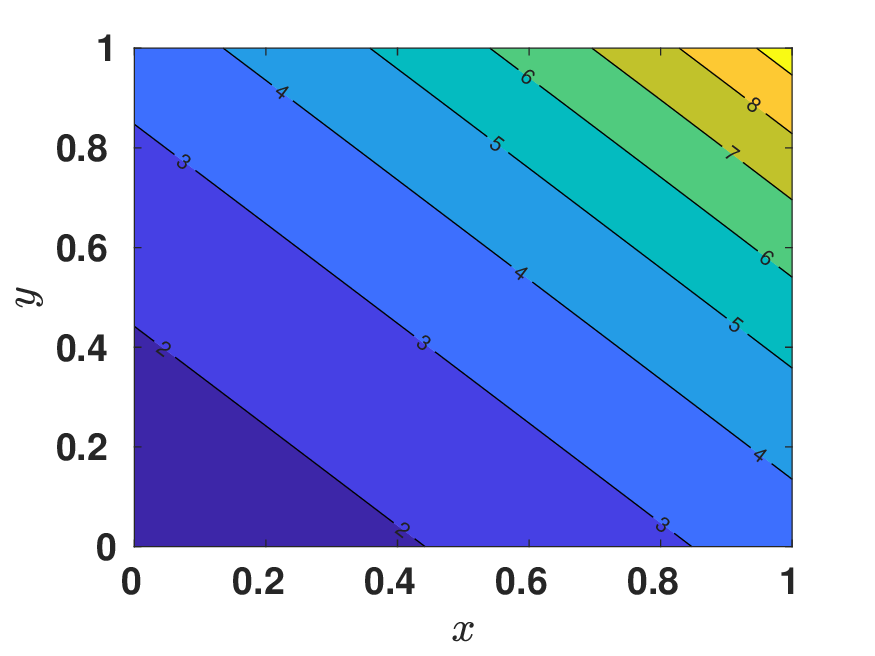}  
		\caption{$t=0.25$}
		\label{fig:Soln_t0.25}
	\end{subfigure}
	\begin{subfigure}{.45\textwidth}
		\centering
		% include second image
		\includegraphics[width=1\linewidth]{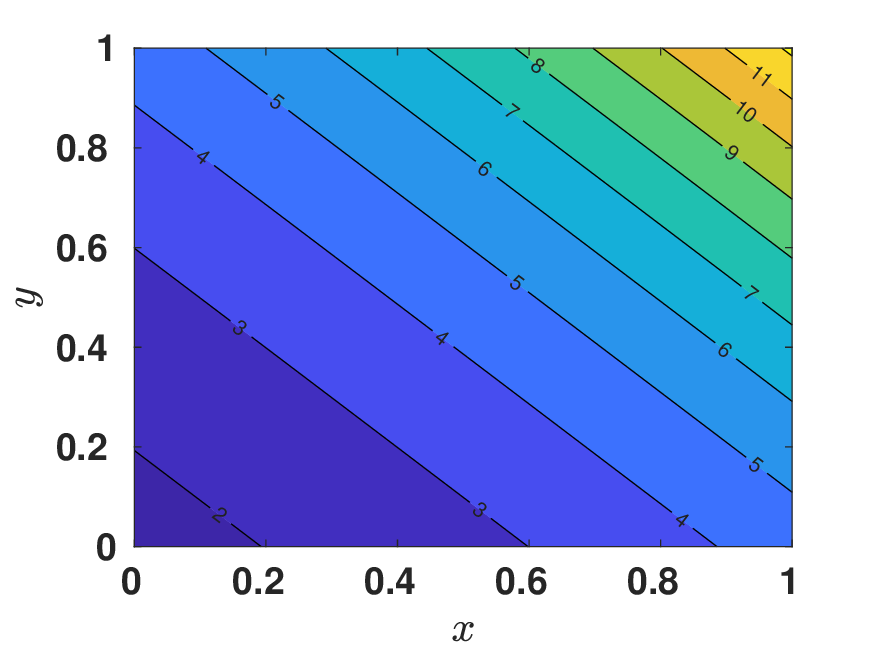}  
		\caption{$t=0.5$}
		\label{fig:Soln_t0.5}
	\end{subfigure}
	
	%		\newline
	
	\begin{subfigure}{.45\textwidth}
		\centering
		% include third image
		\includegraphics[width=1\linewidth]{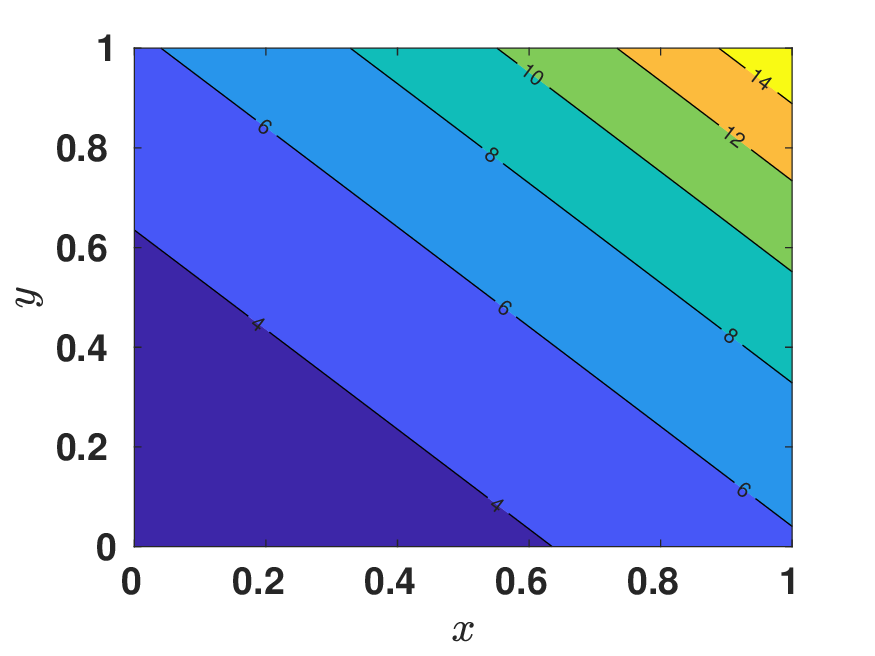}  
		\caption{$t=0.75$}
		\label{fig:Soln_t0.75}
	\end{subfigure}
	\begin{subfigure}{.45\textwidth}
		\centering
		% include fourth image
		\includegraphics[width=1\linewidth]{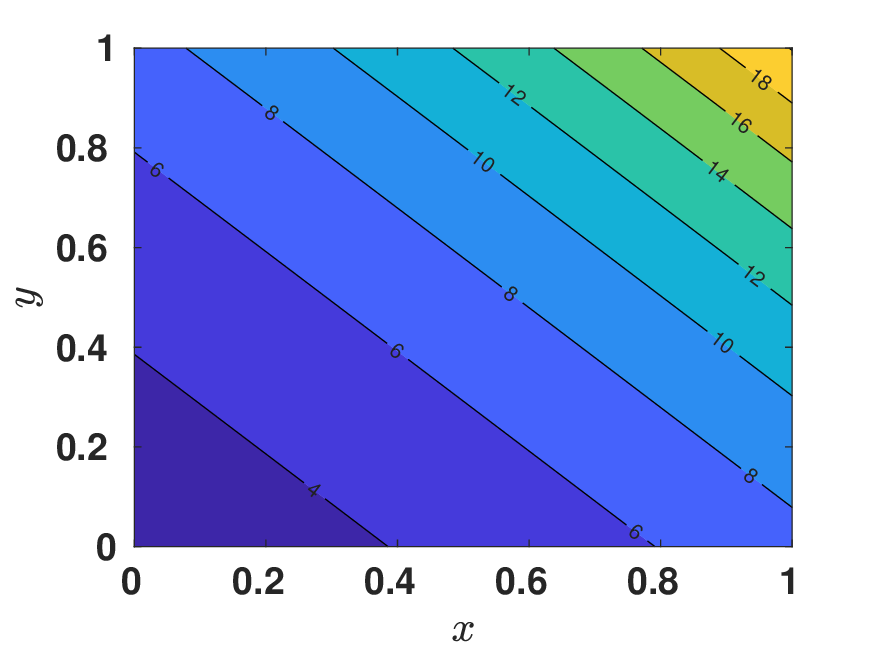}  
		\caption{$t=1$}
		\label{fig:Soln_t1}
	\end{subfigure}
	\caption{Patch dynamics solutions of the problem \eqref{eqn:Physical_Problem} in the domain $[0,1]\times[0,1]$ at various time instances for $\operatorname{N_x}=\operatorname{N_y}=10$, $\operatorname{N_t}=1000$.}
	\label{fig:Soln_Contour}
\end{figure}

Figure \ref{fig:Soln_Contour} shows the snapshots of the patch dynamics solutions of the problem \eqref{eqn:Physical_Problem} in the domain $[0,1]\times[0,1]$ at time instances $t= 0.25$, 0.5, 0.75 and 1. To find the solution, we set $\operatorname{N_x}=\operatorname{N_y}=10$, $\operatorname{N_t}=1000$. The solution shows a high gradient near the top and right boundaries. As time increases, the gradients near the top and right boundaries also increase.

In Table \ref{table:HBMvsPD}, we compare the percentage errors of the problem \eqref{eqn:Physical_Problem} obtained using the $\mathsf{HBM}$ and patch dynamics scheme. In the last four columns, the values represent the percentage error corresponding to each scheme. The results show that the patch dynamics scheme has significantly better accuracy than the $\mathsf{HBM}$ scheme. The percentage errors in the patch dynamics solution are negligible compared to the $\mathsf{HBM}$ solution. For example, the error at the location (0.2,0.2) in the patch dynamics solution is $\frac{1}{3000}$ times the error in the $\mathsf{HBM}$ solution at the exact location, for the same macroscopic spatial grid resolution 20$\times$20.  
This shows that the patch dynamics scheme stands as a powerful tool to solve high gradient, convection-dominated $\mathsf{CDR}$ equations. As the problem \eqref{eqn:Physical_Problem} with \eqref{eqn:P1_P} has a high gradient near the right and top walls (i.e., near $x=1$ and $y=1$) of the domain, the solution contains larger errors in such regions (near (0.4,0.4), (0.6,0.6) and (0.8,0.8)) compared to that in the low gradient regions (near (0.2,0.2)). Similar physical behaviour may also be observed in the boundary layer problems \cite{2001_ge_high}. Such phenomena are more efficiently addressed using appropriate multiscale techniques in space. 

\textbf{$\bullet$ Patch dynamics ($\mathsf{PD}$) solution with transformation:}

To solve the heterogeneous problem \eqref{eqn:Physical_Problem} more accurately, we employ a non-uniform grid as a computational tool to resolve the multiscale features near the top and right boundaries of the physical domain $\Omega_p$. The following orthogonal transformation
\begin{equation}\label{eqn:Stretching_Transformation}
	x=\xi+\frac{\lambda}{\pi}\sin(\pi \xi), \hspace{0.25cm} y=\eta+\frac{\lambda}{\pi}\sin(\pi \eta),
\end{equation}
is used here, where $\lambda$ is the stretching parameter ($0\leq\lambda<1$), which controls the degree of clustering in the domain. $\lambda=0$ shows the uniform grids. The computational domain $[0, 1]\times[0, 1]$ in the $\xi\eta$-plane is converted into the physical domain $[0, 1]\times[0, 1]$ in $xy$-plane under the transformation \eqref{eqn:Stretching_Transformation}. Here, 
\begin{equation}
	J=\left(1+\lambda\cos(\pi\xi)\right)\left(1+\lambda\cos(\pi\eta)\right)\ge (1-\lambda)^2>0,
\end{equation}
throughout the entire domain. The transformation \eqref{eqn:Stretching_Transformation} is smooth, one-to-one, and $0<J<\infty$, which shows that the mapping is non-singular.

\textbf{Computational Problem:}

Under the transformation \eqref{eqn:Stretching_Transformation}, the physical problem \eqref{eqn:Physical_Problem} together with \eqref{eqn:P1_P} is reduced to the computational problem \eqref{eqn:Computational_Problem}, where
\begin{equation}\label{eqn:P1_C}
	\begin{split}
		&\alpha(\xi,\eta)=\frac{1}{\{1+\lambda\cos(\pi\xi)\}^2},\hspace{0.2cm} \beta(\xi,\eta)=0,\hspace{0.2cm} \gamma(\xi,\eta)=\frac{1}{\{1+\lambda\cos(\pi\eta)\}^2}\\
		&\nu(\xi,\eta)=10\frac{\xi+\frac{\lambda}{\pi}\sin(\pi\xi)}{1+\lambda\cos(\pi\xi)}-\frac{\lambda\pi\sin(\pi\xi)}{\{1+\lambda\cos(\pi\xi)\}^3},\\ &\omega(\xi,\eta)=10\frac{\eta+\frac{\lambda}{\pi}\sin(\pi\eta)}{1+\lambda\cos(\pi\eta)}-\frac{\lambda\pi\sin(\pi\eta)}{\{1+\lambda\cos(\pi\eta)\}^3},\\
		&\phi(\xi,\eta)=0, \hspace{0.2cm} g(\xi,\eta,t)= \left[10\left\{\xi+\frac{\lambda}{\pi}\sin(\pi \xi)\right\} 
		+ 10\!\left\{\eta+\frac{\lambda}{\pi}\sin(\pi \eta)\right\} - 1\right] \\
		&\hspace{3cm}\quad \times \exp\!\left[\left\{\xi+\tfrac{\lambda}{\pi}\sin(\pi \xi)\right\} 
		+ \left\{\eta+\tfrac{\lambda}{\pi}\sin(\pi \eta)\right\} + t\right].
	\end{split}
\end{equation}
and the corresponding initial and boundary conditions in $\Omega_c$ are:
\begin{equation}\label{eqn:CDR_Com_IC_BCs}
	\begin{split}
		%		\mathsf{PDE:}\hspace{0.15cm} \frac{\partial u}{\partial t}+v_\xi \frac{\partial u}{\partial \xi}+v_\eta \frac{\partial u}{\partial \eta}=D_\xi\frac{\partial^2 u}{\partial \xi^2}+D_\eta\frac{\partial^2 u}{\partial \eta^2}+S_c(\xi,\eta,t),\hspace{0.25cm} \text{in} \hspace{0.25cm}  \Omega_c = (0,1)\times(0,1),\\
		%		 \hspace{9.6cm}t\in(0,1],\\
		\mathsf{IC:}&\hspace{0.1cm}u(\xi,\eta,0)=\exp\left[\left\{\xi+\frac{\lambda}{\pi}\sin(\pi \xi)\right\}+\left\{\eta+\frac{\lambda}{\pi}\sin(\pi \eta)\right\}\right],\hspace{0.25cm} \text{in}\hspace{0.25cm} \Omega_c,\\
		\mathsf{BCs:}&\hspace{0.1cm}u(0,\eta,t)=\exp\left[\left\{\eta+\frac{\lambda}{\pi}\sin(\pi \eta)\right\}+t\right], \hspace{0.1cm} u(1,y,t)=\exp\left[1+\left\{\eta+\frac{\lambda}{\pi}\sin(\pi \eta)\right\}+t\right],\\
		&u(\xi,0,t)=\exp\left[\left\{\xi+\frac{\lambda}{\pi}\sin(\pi \xi)\right\}+t\right], \hspace{0.1cm} u(\xi,1,t)=\exp\left[\left\{\xi+\frac{\lambda}{\pi}\sin(\pi \xi)\right\}+1+t\right],\\
		&\text{on}\hspace{0.25cm}\partial\Omega_c \hspace{0.25cm} \text{and} \hspace{0.25cm}t \in[0,1].
	\end{split}
\end{equation}
%where,
%\begin{center}
%	$v_\xi=10\frac{\xi+\frac{\lambda}{\pi}\sin(\pi\xi)}{1+\lambda\cos(\pi\xi)}-\frac{\lambda\pi\sin(\pi\xi)}{\{1+\lambda\cos(\pi\xi)\}^3}$, $v_\eta=10\frac{\eta+\frac{\lambda}{\pi}\sin(\pi\eta)}{1+\lambda\cos(\pi\eta)}-\frac{\lambda\pi\sin(\pi\eta)}{\{1+\lambda\cos(\pi\eta)\}^3}$,
%\end{center}  
%\begin{center}
%	$D_\xi=\frac{1}{\{1+\lambda\cos(\pi\xi)\}^2}$, $D_\eta=\frac{1}{\{1+\lambda\cos(\pi\eta)\}^2}$
%\end{center}
%and 
%\[
%\begin{aligned}
%	S_c(\xi,\eta,t) &= \left[10\left\{\xi+\frac{\lambda}{\pi}\sin(\pi \xi)\right\} 
%	+ 10\!\left\{\eta+\frac{\lambda}{\pi}\sin(\pi \eta)\right\} - 1\right] \\[6pt]
%	&\quad \times \exp\!\left[\left\{\xi+\tfrac{\lambda}{\pi}\sin(\pi \xi)\right\} 
%	+ \left\{\eta+\tfrac{\lambda}{\pi}\sin(\pi \eta)\right\} + t\right].
%\end{aligned}
%\]

%\end{center}

\begin{figure}\label{eqn:Uniform_Nonuniform_Grids}
	\centering
	\advance\leftskip-3cm
	\advance\rightskip-3cm
	\begin{tikzpicture}[scale=5]
		% Physical Domain
		\draw[black, very thick] (0,0) rectangle (0.98,0.98);
		
		\draw[->,ultra thick] (-0.1,-0.1)--(0.1,-0.1) node[right]{$x$};
		\draw[->,ultra thick] (-0.1,-0.105)--(-0.1,0.1) node[above]{$y$};
		
		\draw (0,0.1130) -- (0.98,0.1130);
		\draw (0,0.2240) -- (0.98,0.2240);
		\draw (0,0.3310) -- (0.98,0.3310);
		\draw (0,0.4323) -- (0.98,0.4323);
		\draw (0,0.5263) -- (0.98,0.5263);
		\draw (0,0.6119) -- (0.98,0.6119);
		\draw (0,0.6883) -- (0.98,0.6883);
		\draw (0,0.7549) -- (0.98,0.7549);
		\draw (0,0.8119) -- (0.98,0.8119);
		\draw (0,0.8596) -- (0.98,0.8596);
		\draw (0,0.8989) -- (0.98,0.8989);
		\draw (0,0.9310) -- (0.98,0.9310);
		\draw (0,0.9573) -- (0.98,0.9573);
		\draw (0,0.9797) -- (0.98,0.9797);

		\draw (0.1130,0) -- (0.1130,0.98);
		\draw (0.2240,0) -- (0.2240,0.98);
		\draw (0.3310,0) -- (0.3310,0.98);
		\draw (0.4323,0) -- (0.4323,0.98);
		\draw (0.5263,0) -- (0.5263,0.98);
		\draw (0.6119,0) -- (0.6119,0.98);
		\draw (0.6883,0) -- (0.6883,0.98);
		\draw (0.7549,0) -- (0.7549,0.98);
		\draw (0.8119,0) -- (0.8119,0.98);
		\draw (0.8596,0) -- (0.8596,0.98);
		\draw (0.8989,0) -- (0.8989,0.98);
		\draw (0.9310,0) -- (0.9310,0.98);
		\draw (0.9573,0) -- (0.9573,0.98);
		\draw (0.9797,0) -- (0.9797,0.98);

		\draw[fill=orange] (0.075,0.07) rectangle (0.15,0.15);
		\draw[fill=orange] (0.075,0.19) rectangle (0.15,0.26);
		\draw[fill=orange] (0.075,0.3) rectangle (0.15,0.36);
		\draw[fill=orange] (0.075,0.4075) rectangle (0.15,0.4625);
		\draw[fill=orange] (0.075,0.5) rectangle (0.15,0.55);
		\draw[fill=orange] (0.075,0.59) rectangle (0.15,0.63);
		\draw[fill=orange] (0.075,0.672) rectangle (0.15,0.708);
		\draw[fill=orange] (0.075,0.74) rectangle (0.15,0.775);
		\draw[fill=orange] (0.075,0.799) rectangle (0.15,0.821);
		\draw[fill=orange] (0.075,0.845) rectangle (0.15,0.87);
		\draw[fill=orange] (0.075,0.89) rectangle (0.15,0.91);
		\draw[fill=orange] (0.075,0.92) rectangle (0.15,0.938);
		\draw[fill=orange] (0.075,0.949) rectangle (0.15,0.964);
		
		\draw[fill=yellow] (0.19,0.07) rectangle (0.26,0.15);
		\draw[fill=yellow] (0.19,0.19) rectangle (0.26,0.26);
		\draw[fill=yellow] (0.19,0.3) rectangle (0.26,0.36);
		\draw[fill=yellow] (0.19,0.4075) rectangle (0.26,0.4625);
		\draw[fill=yellow] (0.19,0.5) rectangle (0.26,0.55);
		\draw[fill=yellow] (0.19,0.59) rectangle (0.26,0.63);
		\draw[fill=yellow] (0.19,0.672) rectangle (0.26,0.708);
		\draw[fill=yellow] (0.19,0.74) rectangle (0.26,0.775);
		\draw[fill=yellow] (0.19,0.799) rectangle (0.26,0.821);
		\draw[fill=yellow] (0.19,0.845) rectangle (0.26,0.87);
		\draw[fill=yellow] (0.19,0.89) rectangle (0.26,0.91);
		\draw[fill=yellow] (0.19,0.92) rectangle (0.26,0.938);
		\draw[fill=yellow] (0.19,0.949) rectangle (0.26,0.964);
		
		\draw[fill=blue!50] (0.3,0.07) rectangle (0.36,0.15);
		\draw[fill=blue!50] (0.3,0.19) rectangle (0.36,0.26);
		\draw[fill=blue!50] (0.3,0.3) rectangle (0.36,0.36);
		\draw[fill=blue!50] (0.3,0.4075) rectangle (0.36,0.4625);
		\draw[fill=blue!50] (0.3,0.5) rectangle (0.36,0.55);
		\draw[fill=blue!50] (0.3,0.59) rectangle (0.36,0.63);
		\draw[fill=blue!50] (0.3,0.672) rectangle (0.36,0.708);
		\draw[fill=blue!50] (0.3,0.74) rectangle (0.36,0.775);
		\draw[fill=blue!50] (0.3,0.799) rectangle (0.36,0.821);
		\draw[fill=blue!50] (0.3,0.845) rectangle (0.36,0.87);
		\draw[fill=blue!50] (0.3,0.89) rectangle (0.36,0.91);
		\draw[fill=blue!50] (0.3,0.92) rectangle (0.36,0.938);
		\draw[fill=blue!50] (0.3,0.949) rectangle (0.36,0.964);
		
		\draw[fill=red!70] (0.405,0.07) rectangle (0.46,0.15);
		\draw[fill=red!70] (0.405,0.19) rectangle (0.46,0.26);
		\draw[fill=red!70] (0.405,0.3) rectangle (0.46,0.36);
		\draw[fill=red!70] (0.405,0.4075) rectangle (0.46,0.4625);
		\draw[fill=red!70] (0.405,0.5) rectangle (0.46,0.55);
		\draw[fill=red!70] (0.405,0.59) rectangle (0.46,0.63);
		\draw[fill=red!70] (0.405,0.672) rectangle (0.46,0.708);
		\draw[fill=red!70] (0.405,0.74) rectangle (0.46,0.775);
		\draw[fill=red!70] (0.405,0.799) rectangle (0.46,0.821);
		\draw[fill=red!70] (0.405,0.845) rectangle (0.46,0.87);
		\draw[fill=red!70] (0.405,0.89) rectangle (0.46,0.91);
		\draw[fill=red!70] (0.405,0.92) rectangle (0.46,0.938);
		\draw[fill=red!70] (0.405,0.949) rectangle (0.46,0.964);
		
		\draw[fill=teal!70] (0.501,0.07) rectangle (0.551,0.15);
		\draw[fill=teal!70] (0.501,0.19) rectangle (0.551,0.26);
		\draw[fill=teal!70] (0.501,0.3) rectangle (0.551,0.36);
		\draw[fill=teal!70] (0.501,0.4075) rectangle (0.551,0.4625);
		\draw[fill=teal!70] (0.501,0.5) rectangle (0.551,0.55);
		\draw[fill=teal!70] (0.501,0.59) rectangle (0.551,0.63);
		\draw[fill=teal!70] (0.501,0.672) rectangle (0.551,0.708);
		\draw[fill=teal!70] (0.501,0.74) rectangle (0.551,0.775);
		\draw[fill=teal!70] (0.501,0.799) rectangle (0.551,0.821);
		\draw[fill=teal!70] (0.501,0.845) rectangle (0.551,0.87);
		\draw[fill=teal!70] (0.501,0.89) rectangle (0.551,0.91);
		\draw[fill=teal!70] (0.501,0.92) rectangle (0.551,0.938);
		\draw[fill=teal!70] (0.501,0.949) rectangle (0.551,0.964);
		
		\draw[fill=pink!100] (0.592,0.07) rectangle (0.632,0.15);
		\draw[fill=pink!100] (0.592,0.19) rectangle (0.632,0.26);
		\draw[fill=pink!100] (0.592,0.3) rectangle (0.632,0.36);
		\draw[fill=pink!100] (0.592,0.4075) rectangle (0.632,0.4625);
		\draw[fill=pink!100] (0.592,0.5) rectangle (0.632,0.55);
		\draw[fill=pink!100] (0.592,0.59) rectangle (0.632,0.63);
		\draw[fill=pink!100] (0.592,0.672) rectangle (0.632,0.708);
		\draw[fill=pink!100] (0.592,0.74) rectangle (0.632,0.775);
		\draw[fill=pink!100] (0.592,0.799) rectangle (0.632,0.821);
		\draw[fill=pink!100] (0.592,0.845) rectangle (0.632,0.87);
		\draw[fill=pink!100] (0.592,0.89) rectangle (0.632,0.91);
		\draw[fill=pink!100] (0.592,0.92) rectangle (0.632,0.938);
		\draw[fill=pink!100] (0.592,0.949) rectangle (0.632,0.964);
		
		\draw[fill=cyan!100] (0.67,0.07) rectangle (0.706,0.15);
		\draw[fill=cyan!100] (0.67,0.19) rectangle (0.706,0.26);
		\draw[fill=cyan!100] (0.67,0.3) rectangle (0.706,0.36);
		\draw[fill=cyan!100] (0.67,0.4075) rectangle (0.706,0.4625);
		\draw[fill=cyan!100] (0.67,0.5) rectangle (0.706,0.55);
		\draw[fill=cyan!100] (0.67,0.59) rectangle (0.706,0.63);
		\draw[fill=cyan!100] (0.67,0.672) rectangle (0.706,0.708);
		\draw[fill=cyan!100] (0.67,0.74) rectangle (0.706,0.775);
		\draw[fill=cyan!100] (0.67,0.799) rectangle (0.706,0.821);
		\draw[fill=cyan!100] (0.67,0.845) rectangle (0.706,0.87);
		\draw[fill=cyan!100] (0.67,0.89) rectangle (0.706,0.91);
		\draw[fill=cyan!100] (0.67,0.92) rectangle (0.706,0.938);
		\draw[fill=cyan!100] (0.67,0.949) rectangle (0.706,0.964);
		
		\draw[fill=magenta!60] (0.737,0.07) rectangle (0.772,0.15);
		\draw[fill=magenta!60] (0.737,0.19) rectangle (0.772,0.26);
		\draw[fill=magenta!60] (0.737,0.3) rectangle (0.772,0.36);
		\draw[fill=magenta!60] (0.737,0.4075) rectangle (0.772,0.4625);
		\draw[fill=magenta!60] (0.737,0.5) rectangle (0.772,0.55);
		\draw[fill=magenta!60] (0.737,0.59) rectangle (0.772,0.63);
		\draw[fill=magenta!60] (0.737,0.672) rectangle (0.772,0.708);
		\draw[fill=magenta!60] (0.737,0.74) rectangle (0.772,0.775);
		\draw[fill=magenta!60] (0.737,0.799) rectangle (0.772,0.821);
		\draw[fill=magenta!60] (0.737,0.845) rectangle (0.772,0.87);
		\draw[fill=magenta!60] (0.737,0.89) rectangle (0.772,0.91);
		\draw[fill=magenta!60] (0.737,0.92) rectangle (0.772,0.938);
		\draw[fill=magenta!60] (0.737,0.949) rectangle (0.772,0.964);
		
		\draw[fill=lightgray!60] (0.7995,0.07) rectangle (0.823,0.15);
		\draw[fill=lightgray!60] (0.7995,0.19) rectangle (0.823,0.26);
		\draw[fill=lightgray!60] (0.7995,0.3) rectangle (0.823,0.36);
		\draw[fill=lightgray!60] (0.7995,0.4075) rectangle (0.823,0.4625);
		\draw[fill=lightgray!60] (0.7995,0.5) rectangle (0.823,0.55);
		\draw[fill=lightgray!60] (0.7995,0.59) rectangle (0.823,0.63);
		\draw[fill=lightgray!60] (0.7995,0.672) rectangle (0.823,0.708);
		\draw[fill=lightgray!60] (0.7995,0.74) rectangle (0.823,0.775);
		\draw[fill=lightgray!60] (0.7995,0.799) rectangle (0.823,0.821);
		\draw[fill=lightgray!60] (0.7995,0.845) rectangle (0.823,0.87);
		\draw[fill=lightgray!60] (0.7995,0.89) rectangle (0.823,0.91);
		\draw[fill=lightgray!60] (0.7995,0.92) rectangle (0.823,0.938);
		\draw[fill=lightgray!60] (0.7995,0.949) rectangle (0.823,0.964);
		
		\draw[fill=olive!60] (0.847,0.07) rectangle (0.872,0.15);
		\draw[fill=olive!60] (0.847,0.19) rectangle (0.872,0.26);
		\draw[fill=olive!60] (0.847,0.3) rectangle (0.872,0.36);
		\draw[fill=olive!60] (0.847,0.4075) rectangle (0.872,0.4625);
		\draw[fill=olive!60] (0.847,0.5) rectangle (0.872,0.55);
		\draw[fill=olive!60] (0.847,0.59) rectangle (0.872,0.63);
		\draw[fill=olive!60] (0.847,0.672) rectangle (0.872,0.708);
		\draw[fill=olive!60] (0.847,0.74) rectangle (0.872,0.775);
		\draw[fill=olive!60] (0.847,0.799) rectangle (0.872,0.821);
		\draw[fill=olive!60] (0.847,0.845) rectangle (0.872,0.87);
		\draw[fill=olive!60] (0.847,0.89) rectangle (0.872,0.91);
		\draw[fill=olive!60] (0.847,0.92) rectangle (0.872,0.938);
		\draw[fill=olive!60] (0.847,0.949) rectangle (0.872,0.964);
		
		\draw[fill=brown!60] (0.89,0.07) rectangle (0.908,0.15);
		\draw[fill=brown!60] (0.89,0.19) rectangle (0.908,0.26);
		\draw[fill=brown!60] (0.89,0.3) rectangle (0.908,0.36);
		\draw[fill=brown!60] (0.89,0.4075) rectangle (0.908,0.4625);
		\draw[fill=brown!60] (0.89,0.5) rectangle (0.908,0.55);
		\draw[fill=brown!60] (0.89,0.59) rectangle (0.908,0.63);
		\draw[fill=brown!60] (0.89,0.672) rectangle (0.908,0.708);
		\draw[fill=brown!60] (0.89,0.74) rectangle (0.908,0.775);
		\draw[fill=brown!60] (0.89,0.799) rectangle (0.908,0.821);
		\draw[fill=brown!60] (0.89,0.845) rectangle (0.908,0.87);
		\draw[fill=brown!60] (0.89,0.89) rectangle (0.908,0.91);
		\draw[fill=brown!60] (0.89,0.92) rectangle (0.908,0.938);
		\draw[fill=brown!60] (0.89,0.949) rectangle (0.908,0.964);
		
		\draw[fill=darkgray!60] (0.921,0.07) rectangle (0.939,0.15);
		\draw[fill=darkgray!60] (0.921,0.19) rectangle (0.939,0.26);
		\draw[fill=darkgray!60] (0.921,0.3) rectangle (0.939,0.36);
		\draw[fill=darkgray!60] (0.921,0.4075) rectangle (0.939,0.4625);
		\draw[fill=darkgray!60] (0.921,0.5) rectangle (0.939,0.55);
		\draw[fill=darkgray!60] (0.921,0.59) rectangle (0.939,0.63);
		\draw[fill=darkgray!60] (0.921,0.672) rectangle (0.939,0.708);
		\draw[fill=darkgray!60] (0.921,0.74) rectangle (0.939,0.775);
		\draw[fill=darkgray!60] (0.921,0.799) rectangle (0.939,0.821);
		\draw[fill=darkgray!60] (0.921,0.845) rectangle (0.939,0.87);
		\draw[fill=darkgray!60] (0.921,0.89) rectangle (0.939,0.91);
		\draw[fill=darkgray!60] (0.921,0.92) rectangle (0.939,0.938);
		\draw[fill=darkgray!60] (0.921,0.949) rectangle (0.939,0.964);
		
		\draw[fill=violet!60] (0.95,0.07) rectangle (0.965,0.15);
		\draw[fill=violet!60] (0.95,0.19) rectangle (0.965,0.26);
		\draw[fill=violet!60] (0.95,0.3) rectangle (0.965,0.36);
		\draw[fill=violet!60] (0.95,0.4075) rectangle (0.965,0.4625);
		\draw[fill=violet!60] (0.95,0.5) rectangle (0.965,0.55);
		\draw[fill=violet!60] (0.95,0.59) rectangle (0.965,0.63);
		\draw[fill=violet!60] (0.95,0.672) rectangle (0.965,0.708);
		\draw[fill=violet!60] (0.95,0.74) rectangle (0.965,0.775);
		\draw[fill=violet!60] (0.95,0.799) rectangle (0.965,0.821);
		\draw[fill=violet!60] (0.95,0.845) rectangle (0.965,0.87);
		\draw[fill=violet!60] (0.95,0.89) rectangle (0.965,0.91);
		\draw[fill=violet!60] (0.95,0.92) rectangle (0.965,0.938);
		\draw[fill=violet!60] (0.95,0.949) rectangle (0.965,0.964);
		
		% Computational domain
		\draw[step=0.0667cm,very thin] (1.533,0) grid (2.47,1.0005);
		
		\draw[->,ultra thick] (1.433,-0.1)--(1.633,-0.1) node[right]{$\xi$};
		\draw[->,ultra thick] (1.433,-0.1)--(1.433,0.1) node[above]{$\eta$};
		
		\draw[fill=orange] (1.585,0.0517) rectangle (1.615,0.0817);
		\draw[fill=orange] (1.585,0.1183) rectangle (1.615,0.1483);
		\draw[fill=orange] (1.585,0.1850) rectangle (1.615,0.2150);
		\draw[fill=orange] (1.585,0.2517) rectangle (1.615,0.2817);
		\draw[fill=orange] (1.585,0.3183) rectangle (1.615,0.3483);
		\draw[fill=orange] (1.585,0.3850) rectangle (1.615,0.4150);
		\draw[fill=orange] (1.585,0.4517) rectangle (1.615,0.4817);
		\draw[fill=orange] (1.585,0.5183) rectangle (1.615,0.5483);
		\draw[fill=orange] (1.585,0.5850) rectangle (1.615,0.6150);
		\draw[fill=orange] (1.585,0.6517) rectangle (1.615,0.6817);
		\draw[fill=orange] (1.585,0.7183) rectangle (1.615,0.7483);
		\draw[fill=orange] (1.585,0.7850) rectangle (1.615,0.8150);
		\draw[fill=orange] (1.585,0.8517) rectangle (1.615,0.8817);
		\draw[fill=orange] (1.585,0.9183) rectangle (1.615,0.9483);
		
		\draw[fill=yellow] (1.652,0.0517) rectangle (1.682,0.0817);
		\draw[fill=yellow] (1.652,0.1183) rectangle (1.682,0.1483);
		\draw[fill=yellow] (1.652,0.1850) rectangle (1.682,0.2150);
		\draw[fill=yellow] (1.652,0.2517) rectangle (1.682,0.2817);
		\draw[fill=yellow] (1.652,0.3183) rectangle (1.682,0.3483);
		\draw[fill=yellow] (1.652,0.3850) rectangle (1.682,0.4150);
		\draw[fill=yellow] (1.652,0.4517) rectangle (1.682,0.4817);
		\draw[fill=yellow] (1.652,0.5183) rectangle (1.682,0.5483);
		\draw[fill=yellow] (1.652,0.5850) rectangle (1.682,0.6150);
		\draw[fill=yellow] (1.652,0.6517) rectangle (1.682,0.6817);
		\draw[fill=yellow] (1.652,0.7183) rectangle (1.682,0.7483);
		\draw[fill=yellow] (1.652,0.7850) rectangle (1.682,0.8150);
		\draw[fill=yellow] (1.652,0.8517) rectangle (1.682,0.8817);
		\draw[fill=yellow] (1.652,0.9183) rectangle (1.682,0.9483);
		
		\draw[fill=blue!50] (1.72,0.0517) rectangle (1.75,0.0817);
		\draw[fill=blue!50] (1.72,0.1183) rectangle (1.75,0.1483);
		\draw[fill=blue!50] (1.72,0.1850) rectangle (1.75,0.2150);
		\draw[fill=blue!50] (1.72,0.2517) rectangle (1.75,0.2817);
		\draw[fill=blue!50] (1.72,0.3183) rectangle (1.75,0.3483);
		\draw[fill=blue!50] (1.72,0.3850) rectangle (1.75,0.4150);
		\draw[fill=blue!50] (1.72,0.4517) rectangle (1.75,0.4817);
		\draw[fill=blue!50] (1.72,0.5183) rectangle (1.75,0.5483);
		\draw[fill=blue!50] (1.72,0.5850) rectangle (1.75,0.6150);
		\draw[fill=blue!50] (1.72,0.6517) rectangle (1.75,0.6817);
		\draw[fill=blue!50] (1.72,0.7183) rectangle (1.75,0.7483);
		\draw[fill=blue!50] (1.72,0.7850) rectangle (1.75,0.8150);
		\draw[fill=blue!50] (1.72,0.8517) rectangle (1.75,0.8817);
		\draw[fill=blue!50] (1.72,0.9183) rectangle (1.75,0.9483);
		
		\draw[fill=red!70] (1.785,0.0517) rectangle (1.815,0.0817);
		\draw[fill=red!70] (1.785,0.1183) rectangle (1.815,0.1483);
		\draw[fill=red!70] (1.785,0.1850) rectangle (1.815,0.2150);
		\draw[fill=red!70] (1.785,0.2517) rectangle (1.815,0.2817);
		\draw[fill=red!70] (1.785,0.3183) rectangle (1.815,0.3483);
		\draw[fill=red!70] (1.785,0.3850) rectangle (1.815,0.4150);
		\draw[fill=red!70] (1.785,0.4517) rectangle (1.815,0.4817);
		\draw[fill=red!70] (1.785,0.5183) rectangle (1.815,0.5483);
		\draw[fill=red!70] (1.785,0.5850) rectangle (1.815,0.6150);
		\draw[fill=red!70] (1.785,0.6517) rectangle (1.815,0.6817);
		\draw[fill=red!70] (1.785,0.7183) rectangle (1.815,0.7483);
		\draw[fill=red!70] (1.785,0.7850) rectangle (1.815,0.8150);
		\draw[fill=red!70] (1.785,0.8517) rectangle (1.815,0.8817);
		\draw[fill=red!70] (1.785,0.9183) rectangle (1.815,0.9483);
		
		\draw[fill=teal!70] (1.853,0.0517) rectangle (1.883,0.0817);
		\draw[fill=teal!70] (1.853,0.1183) rectangle (1.883,0.1483);
		\draw[fill=teal!70] (1.853,0.1850) rectangle (1.883,0.2150);
		\draw[fill=teal!70] (1.853,0.2517) rectangle (1.883,0.2817);
		\draw[fill=teal!70] (1.853,0.3183) rectangle (1.883,0.3483);
		\draw[fill=teal!70] (1.853,0.3850) rectangle (1.883,0.4150);
		\draw[fill=teal!70] (1.853,0.4517) rectangle (1.883,0.4817);
		\draw[fill=teal!70] (1.853,0.5183) rectangle (1.883,0.5483);
		\draw[fill=teal!70] (1.853,0.5850) rectangle (1.883,0.6150);
		\draw[fill=teal!70] (1.853,0.6517) rectangle (1.883,0.6817);
		\draw[fill=teal!70] (1.853,0.7183) rectangle (1.883,0.7483);
		\draw[fill=teal!70] (1.853,0.7850) rectangle (1.883,0.8150);
		\draw[fill=teal!70] (1.853,0.8517) rectangle (1.883,0.8817);
		\draw[fill=teal!70] (1.853,0.9183) rectangle (1.883,0.9483);
		
		\draw[fill=pink!100] (1.92,0.0517) rectangle (1.95,0.0817);
		\draw[fill=pink!100] (1.92,0.1183) rectangle (1.95,0.1483);
		\draw[fill=pink!100] (1.92,0.1850) rectangle (1.95,0.2150);
		\draw[fill=pink!100] (1.92,0.2517) rectangle (1.95,0.2817);
		\draw[fill=pink!100] (1.92,0.3183) rectangle (1.95,0.3483);
		\draw[fill=pink!100] (1.92,0.3850) rectangle (1.95,0.4150);
		\draw[fill=pink!100] (1.92,0.4517) rectangle (1.95,0.4817);
		\draw[fill=pink!100] (1.92,0.5183) rectangle (1.95,0.5483);
		\draw[fill=pink!100] (1.92,0.5850) rectangle (1.95,0.6150);
		\draw[fill=pink!100] (1.92,0.6517) rectangle (1.95,0.6817);
		\draw[fill=pink!100] (1.92,0.7183) rectangle (1.95,0.7483);
		\draw[fill=pink!100] (1.92,0.7850) rectangle (1.95,0.8150);
		\draw[fill=pink!100] (1.92,0.8517) rectangle (1.95,0.8817);
		\draw[fill=pink!100] (1.92,0.9183) rectangle (1.95,0.9483);
		
		\draw[fill=cyan!100] (1.985,0.0517) rectangle (2.015,0.0817);
		\draw[fill=cyan!100] (1.985,0.1183) rectangle (2.015,0.1483);
		\draw[fill=cyan!100] (1.985,0.1850) rectangle (2.015,0.2150);
		\draw[fill=cyan!100] (1.985,0.2517) rectangle (2.015,0.2817);
		\draw[fill=cyan!100] (1.985,0.3183) rectangle (2.015,0.3483);
		\draw[fill=cyan!100] (1.985,0.3850) rectangle (2.015,0.4150);
		\draw[fill=cyan!100] (1.985,0.4517) rectangle (2.015,0.4817);
		\draw[fill=cyan!100] (1.985,0.5183) rectangle (2.015,0.5483);
		\draw[fill=cyan!100] (1.985,0.5850) rectangle (2.015,0.6150);
		\draw[fill=cyan!100] (1.985,0.6517) rectangle (2.015,0.6817);
		\draw[fill=cyan!100] (1.985,0.7183) rectangle (2.015,0.7483);
		\draw[fill=cyan!100] (1.985,0.7850) rectangle (2.015,0.8150);
		\draw[fill=cyan!100] (1.985,0.8517) rectangle (2.015,0.8817);
		\draw[fill=cyan!100] (1.985,0.9183) rectangle (2.015,0.9483);
		
		\draw[fill=magenta!60] (2.052,0.0517) rectangle (2.082,0.0817);
		\draw[fill=magenta!60] (2.052,0.1183) rectangle (2.082,0.1483);
		\draw[fill=magenta!60] (2.052,0.1850) rectangle (2.082,0.2150);
		\draw[fill=magenta!60] (2.052,0.2517) rectangle (2.082,0.2817);
		\draw[fill=magenta!60] (2.052,0.3183) rectangle (2.082,0.3483);
		\draw[fill=magenta!60] (2.052,0.3850) rectangle (2.082,0.4150);
		\draw[fill=magenta!60] (2.052,0.4517) rectangle (2.082,0.4817);
		\draw[fill=magenta!60] (2.052,0.5183) rectangle (2.082,0.5483);
		\draw[fill=magenta!60] (2.052,0.5850) rectangle (2.082,0.6150);
		\draw[fill=magenta!60] (2.052,0.6517) rectangle (2.082,0.6817);
		\draw[fill=magenta!60] (2.052,0.7183) rectangle (2.082,0.7483);
		\draw[fill=magenta!60] (2.052,0.7850) rectangle (2.082,0.8150);
		\draw[fill=magenta!60] (2.052,0.8517) rectangle (2.082,0.8817);
		\draw[fill=magenta!60] (2.052,0.9183) rectangle (2.082,0.9483);
		
		\draw[fill=lightgray!60] (2.12,0.0517) rectangle (2.15,0.0817);
		\draw[fill=lightgray!60] (2.12,0.1183) rectangle (2.15,0.1483);
		\draw[fill=lightgray!60] (2.12,0.1850) rectangle (2.15,0.2150);
		\draw[fill=lightgray!60] (2.12,0.2517) rectangle (2.15,0.2817);
		\draw[fill=lightgray!60] (2.12,0.3183) rectangle (2.15,0.3483);
		\draw[fill=lightgray!60] (2.12,0.3850) rectangle (2.15,0.4150);
		\draw[fill=lightgray!60] (2.12,0.4517) rectangle (2.15,0.4817);
		\draw[fill=lightgray!60] (2.12,0.5183) rectangle (2.15,0.5483);
		\draw[fill=lightgray!60] (2.12,0.5850) rectangle (2.15,0.6150);
		\draw[fill=lightgray!60] (2.12,0.6517) rectangle (2.15,0.6817);
		\draw[fill=lightgray!60] (2.12,0.7183) rectangle (2.15,0.7483);
		\draw[fill=lightgray!60] (2.12,0.7850) rectangle (2.15,0.8150);
		\draw[fill=lightgray!60] (2.12,0.8517) rectangle (2.15,0.8817);
		\draw[fill=lightgray!60] (2.12,0.9183) rectangle (2.15,0.9483);
		
		\draw[fill=olive!60] (2.186,0.0517) rectangle (2.216,0.0817);
		\draw[fill=olive!60] (2.186,0.1183) rectangle (2.216,0.1483);
		\draw[fill=olive!60] (2.186,0.1850) rectangle (2.216,0.2150);
		\draw[fill=olive!60] (2.186,0.2517) rectangle (2.216,0.2817);
		\draw[fill=olive!60] (2.186,0.3183) rectangle (2.216,0.3483);
		\draw[fill=olive!60] (2.186,0.3850) rectangle (2.216,0.4150);
		\draw[fill=olive!60] (2.186,0.4517) rectangle (2.216,0.4817);
		\draw[fill=olive!60] (2.186,0.5183) rectangle (2.216,0.5483);
		\draw[fill=olive!60] (2.186,0.5850) rectangle (2.216,0.6150);
		\draw[fill=olive!60] (2.186,0.6517) rectangle (2.216,0.6817);
		\draw[fill=olive!60] (2.186,0.7183) rectangle (2.216,0.7483);
		\draw[fill=olive!60] (2.186,0.7850) rectangle (2.216,0.8150);
		\draw[fill=olive!60] (2.186,0.8517) rectangle (2.216,0.8817);
		\draw[fill=olive!60] (2.186,0.9183) rectangle (2.216,0.9483);
		
		\draw[fill=brown!60] (2.252,0.0517) rectangle (2.282,0.0817);
		\draw[fill=brown!60] (2.252,0.1183) rectangle (2.282,0.1483);
		\draw[fill=brown!60] (2.252,0.1850) rectangle (2.282,0.2150);
		\draw[fill=brown!60] (2.252,0.2517) rectangle (2.282,0.2817);
		\draw[fill=brown!60] (2.252,0.3183) rectangle (2.282,0.3483);
		\draw[fill=brown!60] (2.252,0.3850) rectangle (2.282,0.4150);
		\draw[fill=brown!60] (2.252,0.4517) rectangle (2.282,0.4817);
		\draw[fill=brown!60] (2.252,0.5183) rectangle (2.282,0.5483);
		\draw[fill=brown!60] (2.252,0.5850) rectangle (2.282,0.6150);
		\draw[fill=brown!60] (2.252,0.6517) rectangle (2.282,0.6817);
		\draw[fill=brown!60] (2.252,0.7183) rectangle (2.282,0.7483);
		\draw[fill=brown!60] (2.252,0.7850) rectangle (2.282,0.8150);
		\draw[fill=brown!60] (2.252,0.8517) rectangle (2.282,0.8817);
		\draw[fill=brown!60] (2.252,0.9183) rectangle (2.282,0.9483);
		
		\draw[fill=darkgray!60] (2.318,0.0517) rectangle (2.348,0.0817);
		\draw[fill=darkgray!60] (2.318,0.1183) rectangle (2.348,0.1483);
		\draw[fill=darkgray!60] (2.318,0.1850) rectangle (2.348,0.2150);
		\draw[fill=darkgray!60] (2.318,0.2517) rectangle (2.348,0.2817);
		\draw[fill=darkgray!60] (2.318,0.3183) rectangle (2.348,0.3483);
		\draw[fill=darkgray!60] (2.318,0.3850) rectangle (2.348,0.4150);
		\draw[fill=darkgray!60] (2.318,0.4517) rectangle (2.348,0.4817);
		\draw[fill=darkgray!60] (2.318,0.5183) rectangle (2.348,0.5483);
		\draw[fill=darkgray!60] (2.318,0.5850) rectangle (2.348,0.6150);
		\draw[fill=darkgray!60] (2.318,0.6517) rectangle (2.348,0.6817);
		\draw[fill=darkgray!60] (2.318,0.7183) rectangle (2.348,0.7483);
		\draw[fill=darkgray!60] (2.318,0.7850) rectangle (2.348,0.8150);
		\draw[fill=darkgray!60] (2.318,0.8517) rectangle (2.348,0.8817);
		\draw[fill=darkgray!60] (2.318,0.9183) rectangle (2.348,0.9483);
		
		\draw[fill=violet!60] (2.386,0.0517) rectangle (2.416,0.0817);
		\draw[fill=violet!60] (2.386,0.1183) rectangle (2.416,0.1483);
		\draw[fill=violet!60] (2.386,0.1850) rectangle (2.416,0.2150);
		\draw[fill=violet!60] (2.386,0.2517) rectangle (2.416,0.2817);
		\draw[fill=violet!60] (2.386,0.3183) rectangle (2.416,0.3483);
		\draw[fill=violet!60] (2.386,0.3850) rectangle (2.416,0.4150);
		\draw[fill=violet!60] (2.386,0.4517) rectangle (2.416,0.4817);
		\draw[fill=violet!60] (2.386,0.5183) rectangle (2.416,0.5483);
		\draw[fill=violet!60] (2.386,0.5850) rectangle (2.416,0.6150);
		\draw[fill=violet!60] (2.386,0.6517) rectangle (2.416,0.6817);
		\draw[fill=violet!60] (2.386,0.7183) rectangle (2.416,0.7483);
		\draw[fill=violet!60] (2.386,0.7850) rectangle (2.416,0.8150);
		\draw[fill=violet!60] (2.386,0.8517) rectangle (2.416,0.8817);
		\draw[fill=violet!60] (2.386,0.9183) rectangle (2.416,0.9483);
		
	\end{tikzpicture}
	\caption{Left: Physical domain \& Right: Computational domain} \label{fig:P1_Domain_Discretization}
\end{figure}

%As we discussed in Section 2, under the transformation \eqref{eqn:Stretching_Transformation} of the physical problem \eqref{eqn:Physical_Problem} together with \eqref{eqn:P1_P}, the problem \eqref{eqn:Computational_Problem} together with \eqref{eqn:P1_C}, \eqref{eqn:CDR_Com_IC_BCs} is found in the stretched domain with uniform grids. 
The problem \eqref{eqn:Computational_Problem} is considered as the microscopic problem to find the patch dynamics solution at the macroscopic level. A schematic view of the physical domain and the corresponding computational domain is shown in Figure \ref{fig:P1_Domain_Discretization} for $\operatorname{N_x}$=$\operatorname{N_y}$=$\operatorname{N_\xi}$=$\operatorname{N_\eta}$=14. In the physical domain $\Omega_p$, the patches are rectangular in shape, but non-uniform in size. However, under the transformation \eqref{eqn:Stretching_Transformation}, these non-uniform rectangular patches are mapped into uniform rectangular patches in the computational domain $\Omega_c$.

\begin{table}%[h!]
	%	\setstretch{1.2}
	%	\begin{center}
		\caption{Maximum percentage errors of the macroscopic solution for different stretching parameters $\lambda$ and for different resolutions are presented. The resolution is represented in the form $\operatorname{N_x}\times\operatorname{N_y}$, $\operatorname{N_t}$. Here $\operatorname{N_x}, \operatorname{N_y}$ are the number of macro grids along $x$- and $y$-directions in the physical domain $\Omega_p$. $\operatorname{N_t}$ denotes the total number of macro time steps over the entire time interval [0,1].}
		\centering
		\scalebox{1}{
			\begin{tabular}{cccc}
				%				\hline
				%				\multicolumn{1}{|c|}{} & \multicolumn{3}{c|}{Maximum percentage errors for $(\operatorname{N_x}+1)\times(\operatorname{N_y}+1), \operatorname{N_t}$} \\
				\hline
				$\lambda$ & 10$\times$10, 2000 & 15$\times$15, 4500 & 20$\times$20, 8500 \\
				\hline
				0&7.62$\mathrm{e}{-2}$&2.68$\mathrm{e}{-2}$&1.01$\mathrm{e}{-2}$\\
				0.1&2.46$\mathrm{e}{-2}$&9.10$\mathrm{e}{-3}$&7.80$\mathrm{e}{-3}$\\
				0.2&7.74$\mathrm{e}{-2}$&3.78$\mathrm{e}{-2}$&2.42$\mathrm{e}{-2}$\\
				0.3&1.49$\mathrm{e}{-1}$&6.76$\mathrm{e}{-2}$&4.02$\mathrm{e}{-2}$\\
				0.4&2.20$\mathrm{e}{-1}$&9.80$\mathrm{e}{-2}$&5.63$\mathrm{e}{-2}$\\
				0.5&2.91$\mathrm{e}{-1}$&1.29$\mathrm{e}{-1}$&7.27$\mathrm{e}{-2}$\\
				0.6&3.62$\mathrm{e}{-1}$&1.60$\mathrm{e}{-1}$&8.95$\mathrm{e}{-2}$\\
				\hline
		\end{tabular}}
		%	\end{center}
	%\vspace{0.2cm}
	\label{table:Comparison_lambda_grid}
\end{table}

In the computational domain ($\Omega_c$), uniform patch sizes are set as $h_\xi=h_\eta=0.001$, with a time-stepper size of $\tau=1e-6$. At the microscopic level, each patch is discretised using $\operatorname{n}=10$ spatial nano grid steps along both $\xi$ and $\eta$ directions, respectively. The time-stepper is discretised using $\operatorname{n_t}=2$, where $\operatorname{n_t}$ denotes the number of nano time steps. Similarly, a 2-point upwind scheme $[\frac{\partial u}{\partial p}]_{i}=\frac{1}{\delta p}(u_i-u_{i-1})$
of first order and a 3-point central difference scheme of second order are used to discretise the convective and diffusive terms, respectively. $\delta p$ denotes spatial nano step. However, in this convection-dominated problem, higher-order schemes such as the 3-point upwind scheme, given by
$[\frac{\partial u}{\partial p}]_{i}=\frac{1}{2 \hspace{0.05cm}\delta p}(3u_i-4u_{i-1}+u_{i-2})$ discussed in \cite{2020_anderson_computational},
and the 4-point upwind scheme, given by
$[\frac{\partial u}{\partial p}]_{i}=\frac{1}{2\hspace{0.05cm}\delta p}(u_{i+1}-u_{i-1})+\frac{r}{3\hspace{0.05cm}\delta p}(u_{i-2}-3u_{i-1}+3u_{i}-u_{i+1})$ discussed in \cite{2012_fletcher_computational},
both of second-order accurate, were also used.
%Here, the parameter $r$ controls the size of the modification. The 3-point and 4-point upwind schemes exhibit a higher order of accuracy compared to the 2-point upwind scheme. However, 
At macroscopic level, the patch dynamics solution yields comparable accuracy for all three schemes. The 3-point and 4-point upwind schemes incur higher computational costs; so, the 2-point upwind scheme proves to be more efficient. This observation supports the comment made by Maclean et al. \cite{2015_maclean_convergence}, where they noted that increasing the order of the microsolver does not necessarily improve the overall accuracy of the patch dynamics solution. 

Table \ref{table:Comparison_lambda_grid} presents the maximum percentage errors in the physical domain $\Omega_p$ over the full-time interval [0,1] for different spatial resolutions, namely, 10$\times$10, 15$\times$15 and 20$\times$20, under varying values of the clustering parameter $\lambda$. For $0\leq\lambda\leq 0.6$, the results show a good agreement with the exact solution, with the best accuracy consistently achieved at $\lambda=0.1$ across all grid sizes.
This indicates that a clustered grid with $\lambda=0.1$ performs better than a uniform grid (i.e. $\lambda=0$). Furthermore, the results in Table \ref{table:Comparison_lambda_grid} demonstrate that the patch dynamics solutions converge to the exact solution as the macro grid is refined, thereby confirming the convergence of the patch dynamics scheme. From the numerical experiments, we observed that the computational time and the memory requirements of the patch dynamics scheme do not depend noticeably on the stretching parameter $\lambda$.

%%%%%%%%%%%%%%%%%%%%%%%%%%%%%%%%%%%%%%%%%%%%%%%%%%%%%

\begin{figure}%[h!]
	\begin{subfigure}{.5\textwidth}
		\centering
		% include first image
		\includegraphics[width=1\linewidth]{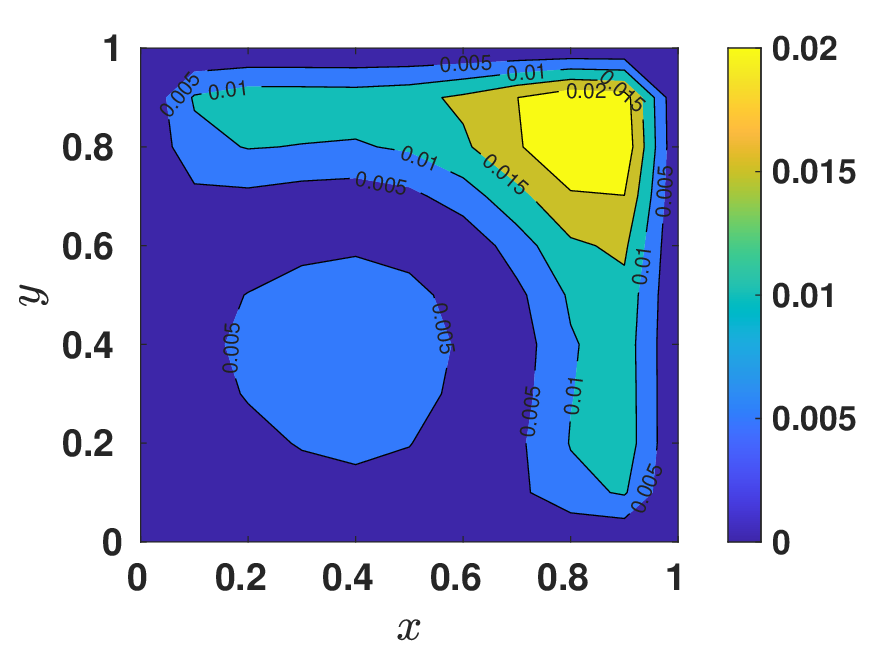}  
		\caption{Percentage error for 10$\times$10 grid.}
		\label{fig:CDR_Contour_11_11}
	\end{subfigure}
	\begin{subfigure}{.5\textwidth}
		\centering
		% include second image
		\includegraphics[width=1\linewidth]{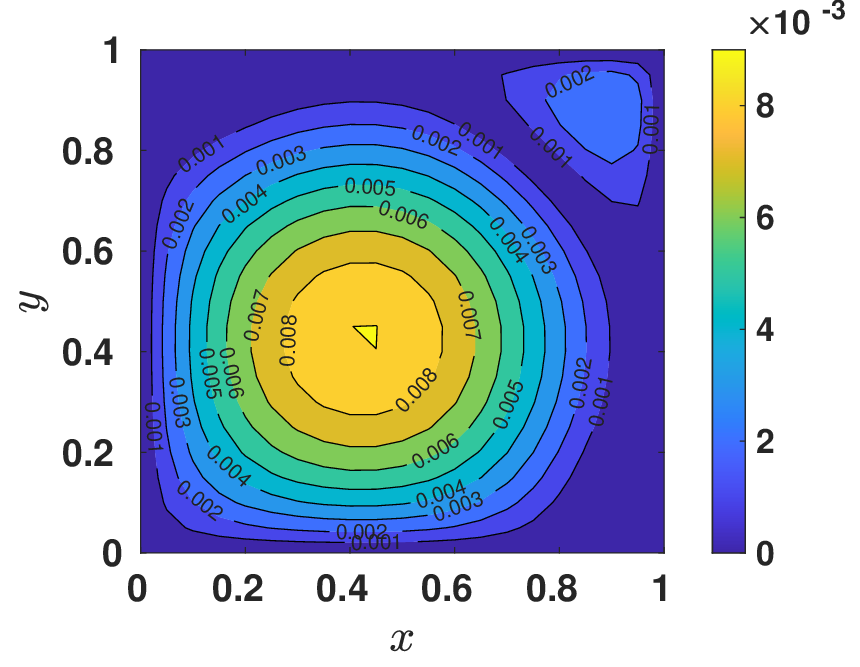}  
		\caption{Percentage error for 20$\times$20 grid.}
		\label{fig:CDR_Contour_21_21}
	\end{subfigure}
	\caption{Contour plots of the percentage errors of the patch dynamics solutions with congregated grid resolutions 10$\times$10, $\operatorname{N_t}=2000$ and 20$\times$20, $\operatorname{N_t}=8500$ in the physical domain at time $t=1$ are presented here, with $\lambda=0.1$ in both.}
	\label{fig:P1_CDR_Contour}
\end{figure}

In Figure \ref{fig:P1_CDR_Contour}, a comparison is presented between the percentage errors of two solutions computed with spatial resolutions of 10$\times$10, $\operatorname{N_t}=2000$ and 20$\times$20, $\operatorname{N_t}=8500$, respectively. In both cases, the stretching parameter is set to $\lambda=0.1$, and all other parameters remain the same, as discussed previously. In Figure \ref{fig:CDR_Contour_11_11}, the solution exhibits excellent overall accuracy, although minor errors are observed in regions with high gradients, which is attributed to the coarser macro grid. In contrast, Figure \ref{fig:CDR_Contour_21_21} demonstrates improved accuracy in the high gradient region due to the use of a finer, more clustered grid. Overall, both solutions demonstrate excellent accuracy.

\begin{table}%[h!]
	%	\setstretch{1.2}
%	\begin{center}
	\caption{%\setstretch{1.2}
		The second column presents the maximum percentage errors in the patch dynamics solution ($U$) relative to the exact solution ($u_e$) for the heterogeneous problem \eqref{eqn:Physical_Problem} with \eqref{eqn:P1_P}, corresponding to various micro-simulation times ($\tau$) listed in the first column.}
	\centering
		\scalebox{1}{
			%	\centering
			
			\begin{tabular}{cc}
				%				\hline
				%				\multicolumn{1}{|c|}{} & \multicolumn{2}{c|}{Maximum percentage errors for $(\operatorname{N_x}+1)\times(\operatorname{N_y}+1), \operatorname{N_t}$} \\
				\hline
				$\tau$ & Max. Per. Error\\ %$\left\|\frac{U-u_e}{u_e}\right\|_\infty\times100\%$ \\
				\hline
				$1\mathrm{e}{-6}$&0.080\\
				$2\mathrm{e}{-6}$&0.085\\
				$4\mathrm{e}{-6}$&0.087\\
				$6\mathrm{e}{-6}$&0.087\\
				$8\mathrm{e}{-6}$&0.088\\
				$10\mathrm{e}{-6}$&0.088\\
				\hline
		\end{tabular}}
%	\end{center}
	%\vspace{0.2cm}
	\label{table:relaxation_time_lambda0}
\end{table}

For problem \eqref{eqn:Physical_Problem} with \eqref{eqn:P1_P}, Table \ref{table:relaxation_time_lambda0} compares the accuracy of the patch dynamics solutions for different micro time duration $\tau$. For this experiment, we set $\operatorname{N_x}=\operatorname{N_y}=10$, $\operatorname{N_t}=2000$, $\lambda=0$, and used same microscale discretisation as described in Sub-section 4.1.1. During this experiment, we varied only the micro time step $\tau$ while keeping all other parameters fixed. As shown in Table \ref{table:relaxation_time_lambda0}, increasing the micro time step ($\tau$) results in no noticeable change in the maximum percentage error. These findings agree with the results reported by Samaey et al. \cite{2005_samaey_gap} in the context of the gap-tooth scheme, where a smaller time step indeed produces a smaller error. Samaey et al. \cite{2006_samaey_patch} also made similar observations in the context of the patch dynamics scheme. Based on the above discussions, we chose the micro time step $\tau=1e-6$ for all subsequent experiments on the heterogeneous $\mathsf{CDR}$ problems in this Sub-section 4.1.

\vspace{0.5cm}

$\bullet$ \textbf{Grid independence of the heterogeneous $\mathsf{CDR}$ $\mathsf{PDE}$ in stretched domain:}

Let $U_{N\times N}$ and $U_{2N\times 2N}$ be the patch dynamics solutions with respect to the $N\times N$ and $2N\times 2N$ macroscopic grids, respectively, at a certain time. Relative difference of the above two solutions is $\left|\frac{U_{2N\times 2N}-U_{N\times N}}{U_{2N\times 2N}}\right|$, at any time $t$.

To check the grid independence of the problem \eqref{eqn:Physical_Problem} together with \eqref{eqn:P1_P} , the coarser grid $10\times10$, $\operatorname{N_t}=2000$ and the finer grid $20\times20$, $\operatorname{N_t}=8500$ are considered. The stretching factor is $\lambda=0.1$, and all other parameters are kept fixed in both solutions as contemplated in the previous. The maximum relative difference between the coarser and finer grid solutions is found to be 2.34$\mathrm{e}{-4}$ over the whole physical domain [0, 1]$\times$[0, 1] and throughout the full-time interval [0, 1]. If the maximum relative difference between the coarser and finer grid solutions is less than 1$\mathrm{e}{-3}$, then it would be declared as grid independent \cite{2015_mazumder_numerical}. So, the patch dynamics solution with the grid resolution $20\times20$, $\operatorname{N_t}=8500$ is considered as grid independent.

\subsubsection{Patch dynamics ($\mathsf{PD}$) solution with variable diffusivity:}
%\textbf{$\mathsf{CDR}$ equation with variable diffusivity:}
A class of materials called Functionally Graded Materials ($\mathsf{FGMs}$) was invented by Niino et al. \cite{1987_niino_functionally} in 1984 that could withstand extreme temperature gradients. In general, $\mathsf{FGMs}$ can be divided into three categories: gradient composition, gradient porosity and gradient microstructure %, where the microstructure can be found in the articles 
\cite{2021_mohammadi_functionally,2017_mahmoud_lattice}. %Looking into the micro-structure of the $\mathsf{FGMs}$, the temperature distribution is solved in a large domain over a long time using micro-simulations in small patches in the following problem.
In order to study the temperature distribution of the $\mathsf{FGMs}$ on a large domain over a long time, the following problem is selected. %This problem considers the temperature distribution informations of the microstructure, which is taken care by the micro-simulations in the patch dynamics scheme. To investigate how the patch dynamics scheme adapts to $\mathsf{CDR}$ problems in functionally graded media, a model with variable velocity and diffusivity is considered. 
In the physical problem \eqref{eqn:Physical_Problem}, we assume that the thermal diffusivity varies with the horizontal distance, which is $D(x,y)=1+x$, convective velocity is $v=\begin{bmatrix}
	1+10x\hspace{0.3cm} 10y
\end{bmatrix}^\top$, the reaction coefficient and the source terms are $f$=20 and $g(x,y,t)=(8x+10y-1)\exp(x+y+t)$, respectively. Equation \eqref{eqn:Physical_Problem} represents a heterogeneous $\mathsf{CDR}$ equation, as both the diffusion coefficient and the convective velocity depend explicitly on the spatial variables and the source term depends on space-time variables. The initial and boundary conditions are considered to be the same as in Sub-section 4.1.1. This physical problem has an exact solution \eqref{eqn:CDR_Phy_Exact}. Based on the physical behaviour of the problem, a clustered grid is employed in the higher-gradient region using the transformation \eqref{eqn:Stretching_Transformation}, which provides a better accuracy compared to uniform (unclustered) grids. 

Under the transformation \eqref{eqn:Stretching_Transformation}, the physical problem \eqref{eqn:Physical_Problem} is reduced to the computational problem \eqref{eqn:Computational_Problem}, where
\begin{equation}\label{eqn:P1_C_Vari_Diff}
	\begin{split}
		&\alpha(\xi,\eta)=\frac{1+\left\{\xi+\frac{\lambda}{\pi}\sin(\pi\xi)\right\}}{\{1+\lambda\cos(\pi\xi)\}^2},\hspace{0.2cm} \beta(\xi,\eta)=0,\hspace{0.2cm} \gamma(\xi,\eta)=\frac{1+\left\{\xi+\frac{\lambda}{\pi}\sin(\pi\xi)\right\}}{\{1+\lambda\cos(\pi\eta)\}^2}\\
		&\nu(\xi,\eta)=10\frac{\xi+\frac{\lambda}{\pi}\sin(\pi\xi)}{1+\lambda\cos(\pi\xi)}-\frac{\lambda\pi\sin(\pi\xi)\left[1+\left\{\xi+\frac{\lambda}{\pi}\sin(\pi\xi)\right\}\right]}{\{1+\lambda\cos(\pi\xi)\}^3},\\ &\omega(\xi,\eta)=10\frac{\eta+\frac{\lambda}{\pi}\sin(\pi\eta)}{1+\lambda\cos(\pi\eta)}-\frac{\lambda\pi\sin(\pi\eta)\left[1+\left\{\xi+\frac{\lambda}{\pi}\sin(\pi\xi)\right\}\right]}{\{1+\lambda\cos(\pi\eta)\}^3},\\
		&\phi(\xi,\eta)=0, \hspace{0.2cm} g(\xi,\eta,t)= \left[8\left\{\xi+\frac{\lambda}{\pi}\sin(\pi \xi)\right\} 
		+ 10\!\left\{\eta+\frac{\lambda}{\pi}\sin(\pi \eta)\right\} - 1\right] \\
		&\hspace{3cm}\quad \times \exp\!\left[\left\{\xi+\tfrac{\lambda}{\pi}\sin(\pi \xi)\right\} 
		+ \left\{\eta+\tfrac{\lambda}{\pi}\sin(\pi \eta)\right\} + t\right].
	\end{split}
\end{equation}  %$\xi+\frac{\lambda}{\pi}\sin(\pi\xi)$ $\eta+\frac{\lambda}{\pi}\sin(\pi\eta)$ $1+\lambda\cos(\pi\xi)$ $1+\lambda\cos(\pi\eta)$
and the corresponding initial and boundary conditions are the same as those in equation \eqref{eqn:CDR_Com_IC_BCs}.
%\begin{equation}\label{eqn:CDR_Com_IC_BCs}
%	\begin{split}
	%		%		\mathsf{PDE:}\hspace{0.15cm} \frac{\partial u}{\partial t}+v_\xi \frac{\partial u}{\partial \xi}+v_\eta \frac{\partial u}{\partial \eta}=D_\xi\frac{\partial^2 u}{\partial \xi^2}+D_\eta\frac{\partial^2 u}{\partial \eta^2}+S_c(\xi,\eta,t),\hspace{0.25cm} \text{in} \hspace{0.25cm}  \Omega_c = (0,1)\times(0,1),\\
	%		%		 \hspace{9.6cm}t\in(0,1],\\
	%		\mathsf{IC:}&\hspace{0.1cm}u(\xi,\eta,0)=\exp\left[\left\{\xi+\frac{\lambda}{\pi}\sin(\pi \xi)\right\}+\left\{\eta+\frac{\lambda}{\pi}\sin(\pi \eta)\right\}\right],\hspace{0.25cm} \text{in}\hspace{0.25cm} \Omega_c,\\
	%		\mathsf{BCs:}&\hspace{0.1cm}u(0,\eta,t)=\exp\left[\left\{\eta+\frac{\lambda}{\pi}\sin(\pi \eta)\right\}+t\right], \hspace{0.1cm} u(1,y,t)=\exp\left[1+\left\{\eta+\frac{\lambda}{\pi}\sin(\pi \eta)\right\}+t\right],\\
	%		&u(\xi,0,t)=\exp\left[\left\{\xi+\frac{\lambda}{\pi}\sin(\pi \xi)\right\}+t\right], \hspace{0.1cm} u(\xi,1,t)=\exp\left[\left\{\xi+\frac{\lambda}{\pi}\sin(\pi \xi)\right\}+1+t\right],\\
	%		&\text{on}\hspace{0.25cm}\partial\Omega_c \hspace{0.25cm} \text{and} \hspace{0.25cm}t \in[0,1].
	%	\end{split}
%\end{equation}
A schematic representation of the physical and computational domains, along with patch arrangement, is shown in Figure \ref{fig:P1_Domain_Discretization}. A similar multiscale behaviour in both space and time is observed in the $\mathsf{CDR}$ equation with variable diffusivity, as discussed for the problem \eqref{eqn:Physical_Problem} with constant diffusivity in Sub-section 4.1.1.

\begin{table}%[h!]
%		\begin{center}
	\caption{%\setstretch{1.2}
		Maximum percentage errors in the patch dynamics solution for different stretching parameters $\lambda$ and for different resolutions are presented. The resolution is represented in the form of $\operatorname{N_x}\times\operatorname{N_y}$, $\operatorname{N_t}$. Here $\operatorname{N_x}, \operatorname{N_y}$ are the number of macro grids along $x$- and $y$-directions in the physical domain $\Omega_p$. $\operatorname{N_t}$ denotes the total number of macro time steps over the entire interval [0, 1].}
	\centering
		\scalebox{1}{
			%	\centering
			
			\begin{tabular}{ccc}
				%				\hline
				%				& \multicolumn{2}{c|}{Maximum percentage errors} \\
				\hline
				$\lambda$ & 10$\times$10, 2000 & 15$\times$15, 5500\\
				\hline
				0&5.30$\mathrm{e}{-2}$&1.85$\mathrm{e}{-2}$\\
				0.1&1.02$\mathrm{e}{-2}$&8.60$\mathrm{e}{-3}$\\
				0.2&6.67$\mathrm{e}{-2}$&3.22$\mathrm{e}{-2}$\\
				0.3&1.22$\mathrm{e}{-1}$&5.53$\mathrm{e}{-2}$\\
				0.4&1.74$\mathrm{e}{-1}$&7.79$\mathrm{e}{-2}$\\
				0.5&2.24$\mathrm{e}{-1}$&1.00$\mathrm{e}{-1}$\\
				\hline
		\end{tabular}}
%			\end{center}
	%\vspace{0.2cm}
	\label{table:Comparison_lambda_grid_Vari_Diff}
\end{table}

%\begin{table}%[h!]
%	%	\setstretch{1.2}
%	%	\begin{center}
	%		\scalebox{1}{
		%			%	\centering
		%			
		%			\begin{tabular}{|p{0.5cm}|p{1cm}|p{1cm}|p{1cm}|p{1cm}|p{1cm}|p{1cm}|}
			%%				\hline
			%%				$\lambda$ & 10$\times$10, 2000 & 15$\times$15, 5500\\
			%        	    \hline
			%        	    \multicolumn{1}{|c|}{} &\multicolumn{3}{c|}{10$\times$10, 2000} & \multicolumn{3}{c|}{15$\times$15, 5500} \\
			%				\hline
			%				$\lambda$ &Max \%err&Time ($s$)&Memory (MB) & Max \%err&Time&Memory\\
			%				\hline
			%				0&5.3$\mathrm{e}{-2}$&196&3.7&1.9$\mathrm{e}{-2}$&&\\
			%				0.1&1.0$\mathrm{e}{-2}$&198&3.7&8.6$\mathrm{e}{-3}$&&\\
			%				0.2&6.7$\mathrm{e}{-2}$&200&3.7&3.2$\mathrm{e}{-2}$&&\\
			%				0.3&1.2$\mathrm{e}{-1}$&200&3.7&5.5$\mathrm{e}{-2}$&&\\
			%				0.4&1.7$\mathrm{e}{-1}$&200&3.7&7.8$\mathrm{e}{-2}$&&\\
			%				0.5&2.2$\mathrm{e}{-1}$&200&3.7&1.0$\mathrm{e}{-1}$&&\\
			%				\hline
			%		\end{tabular}}
	%		%	\end{center}
%	%\vspace{0.2cm}
%	\caption{%\setstretch{1.2}
	%		Maximum percentage errors in the patch dynamics solution for different stretching ratios $\lambda$ and for different resolutions are presented. The resolution is represented in the form of $\operatorname{N_x}\times\operatorname{N_y}$, $\operatorname{N_t}$. Here $\operatorname{N_x}, \operatorname{N_y}$ are the number of macro grids along $x$ and $y$ direction in the physical domain $\Omega_p$. $\operatorname{N_t}$ denotes the total number of macro time steps over the entire interval [0, 1].}
%	\label{table:Comparison_lambda_grid_Vari_Diff}
%\end{table}	

\begin{figure}%[h!]
	\begin{subfigure}{.5\textwidth}
		\centering
		% include first image
		\includegraphics[width=1\linewidth]{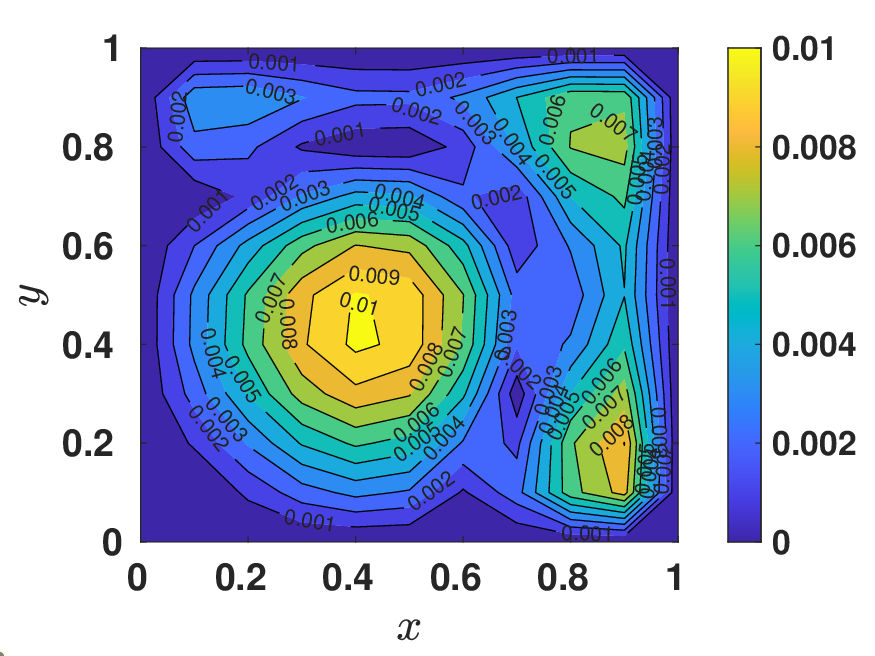}  
		\caption{Percentage error for 10$\times$10 grid.}
		\label{fig:CDR_Contour_VariDiff_11_11}
	\end{subfigure}
	\begin{subfigure}{.5\textwidth}
		\centering
		% include second image
		\includegraphics[width=1\linewidth]{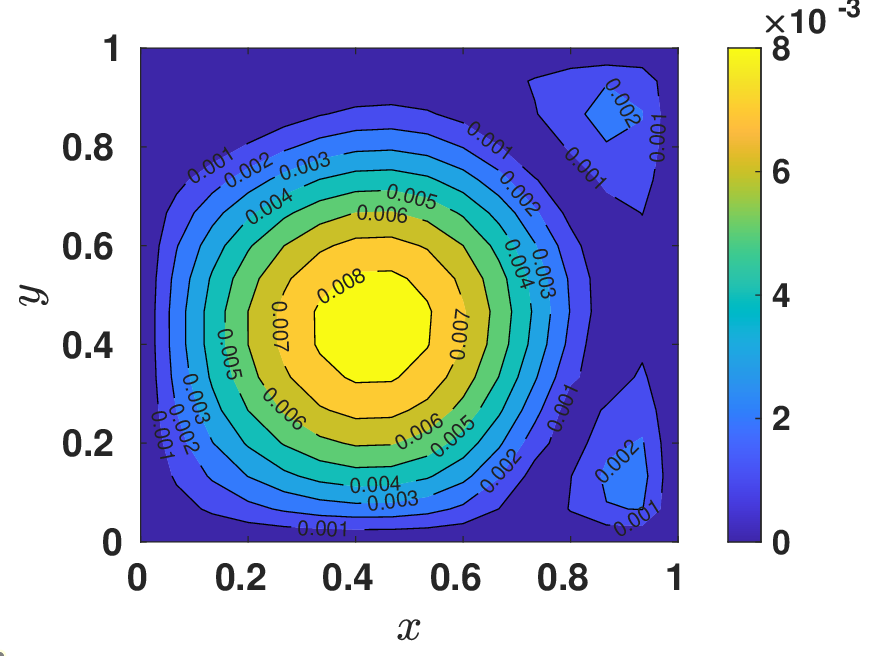}  
		\caption{Percentage error for 20$\times$20 grid.}
		\label{fig:CDR_Contour_VariDiff_21_21}
	\end{subfigure}
	\caption{Contour plots of the percentage errors of the patch dynamics solutions with congregated grids 10$\times$10, $\operatorname{N_t}=2000$ and 20$\times$20, $\operatorname{N_t}=8500$ in the physical domain at time $t=1$ are shown here, with $\lambda=0.1$ in both of the solutions.}
	\label{fig:P1_CDR_Vari_Contour}
\end{figure}

In the variable diffusivity problem, Table \ref{table:Comparison_lambda_grid_Vari_Diff} shows that the best possible solution is obtained for the stretched grid with stretching parameter $\lambda=0.1$. In Figure \ref{fig:P1_CDR_Vari_Contour}, contour plots of the percentage errors of the patch dynamics solutions for congregated grids $10\times10$ and $20\times20$ are shown at the final time $t=1$ and with $\lambda=0.1$ in both the solutions.

%%%%%%%%%%%%%%%%%%%%%%%%%%%%%%%%%%%%%%%%%%%%%%%%%%%%%

$\bullet$ \textbf{Grid independence of the heterogeneous $\mathsf{CDR}$ $\mathsf{PDE}$ in stretched domain:}

Grid independence is examined for the variable diffusivity problem using a coarser grid of $10\times10$, $\operatorname{N_t}=2000$, and a finer grid of $20\times20$, $\operatorname{N_t}=8500$. In both cases, the stretching factor is set to $\lambda=0.1$, and all other parameters are kept the same as previously described. The maximum relative difference between the coarser and finer grid solutions is found to be 9.3$\mathrm{e}{-5}$ (i.e., less than 1$\mathrm{e}{-3}$) over the entire physical domain $[0, 1]^2$ and throughout the full-time interval [0, 1]. Since the maximum relative difference is below 1$\mathrm{e}{-3}$, the patch dynamics solution with the grid resolution $20\times20$ and $\operatorname{N_t}=8500$ is considered to be grid independent \cite{2015_mazumder_numerical}.

$\bullet$ \textbf{A comparison between the performances of the patch dynamics scheme and the full-domain simulation:}

In order to perform this experiment, we consider the heterogeneous problem \eqref{eqn:Physical_Problem} with coefficients as described in the Subsection 4.1.2, where heterogeneity is present in diffusion, advection and source terms. We set $\operatorname{N_x}=\operatorname{N_y}=10$, $\operatorname{N_t}=2000$, $h_\xi=h_\eta=1\mathrm{e}{-3}$, $\tau=1\mathrm{e}{-6}$ and $\operatorname{n}=10$ to compute the patch dynamics solution. A full-domain microsimulation with micro-steps $\delta x=\delta y=1\mathrm{e}{-4}$ is computationally infeasible on a standard computer, as it requires $(10^4+1)^2$ lattice points. Using classical numerical approaches, our computer equipped with 15.7 $\mathsf{GB}$ of $\mathsf{RAM}$ is capable of handling up to $\operatorname{n}=150$, i.e., $151^2=22801$ lattice points over the full domain $[0,1]^2$. Similarly, the computer is unable to simulate the problem using a micro time step of $\delta t=5\mathrm{e}{-7}$. To address these limitations, we employed the method of lines to solve the heterogeneous problem over the entire domain $[0,1]^2$ and full time interval $[0,1]$. This approach is commonly used in equation-free frameworks \cite{2020_Roberts_toolbox,2026_KUMARKARMAKAR_GPD}. The diffusion and advection terms are discretised using a 3-point central difference scheme and a 2-point upwind scheme respectively, which are consistent with the patch dynamics discretisation as described in the Subsection 4.1.1. The resulting semi-discretised system of $\mathsf{ODEs}$ are solved using a variable order and variable step time integration method $\mathsf{ode15s}$. 

\begin{table}%[h!]
	%	\setstretch{1.2}
%	\begin{center}
	\caption{%\setstretch{1.2}
		The performances of the patch dynamics ($\mathsf{PD}$) scheme and the full-domain simulation ($\mathsf{FDS}$) are compared for the heterogeneous problem (1), in which the heterogeneity appears in the diffusion, convective velocity, and source terms. The comparison is carried out in terms of accuracy, computational time and memory usage. The notations $\lambda$, $\%E$, $T$ and $M$ denote the stretching parameter, maximum percentage error over the entire spatial domain $[0,1]^2$ \& time interval $[0,1]$, computational time and the memory requirement, respectively. }
	\centering
		\scalebox{1}{
			%	\centering
			
			\begin{tabular}{cccc}
				%				\hline
				%				\multicolumn{1}{|c|}{} & \multicolumn{2}{c|}{Maximum percentage errors for $(\operatorname{N_x}+1)\times(\operatorname{N_y}+1), \operatorname{N_t}$} \\
				\hline
				$\lambda$ & $\%E_\mathsf{FDS}/\%E_\mathsf{PD}$ & $T_\mathsf{FDS}/T_\mathsf{PD}$ & $M_\mathsf{FDS}/M_\mathsf{PD}$  \\
				\hline
				0&2.7&1.8&1.1\\
				0.1&12.7&2.5&1.1\\
				0.2&1.7&2.6&1.2\\
				\hline
		\end{tabular}}
%	\end{center}
	%\vspace{0.2cm}
	\label{table:PDvsFDS}
\end{table}

Table \ref{table:PDvsFDS} presents a comparison between the performances of the patch dynamics ($\mathsf{PD}$) scheme and the full-domain simulation ($\mathsf{FDS}$). The comparison is based on the accuracy ($\%E$), measured in terms of the maximum percentage error, as well as computational time ($T$) and memory usage ($M$), which are reported in the last three columns of the table. The experiments are conducted for the stretching parameter $\lambda=0, 0.1$ and 0.2, as listed in the first row. The results demonstrate that the patch dynamics scheme consistently achieves a better accuracy compared to the full domain simulation. In addition, the patch dynamics scheme requires less computational time and memory. For instance, when $\lambda=0.1$, the patch dynamics solution attains approximately $\frac{1}{13}$ of the maximum percentage error observed in the full-domain simulation. Moreover, it requires about $\frac{2}{5}$ less computational time and approximately $10\%$ less memory compared to the full-domain simulation. In case of spatial discretisation, the patch dynamics scheme employs a total of $(11\times11)\times(9\times9)+40=9841$ macro-micro lattices, whereas the full-domain simulation uses $151\times151=22801$ lattice points over the entire domain $[0,1]^2$. The full-domain simulation utilises a larger number of lattice points and a higher-order time integration scheme, whereas the proposed patch dynamics scheme delivers superior performance while using significantly fewer lattice points and only a first-order time integration scheme at the macroscopic level.

Overall, considering accuracy, computational time and memory usage, the patch dynamics scheme performs much better compared to the full-domain simulation. 

$\bullet$ \textbf{The proposed patch dynamics scheme is applicable to both $\mathsf{type-A}$ and $\mathsf{type-B}$ problems:}

Weinan E. \cite{2003_weinan_heterognous,2011_weinan_principles} classified multiscale problems into two categories: $\mathsf{type-A}$ and $\mathsf{type-B}$. In $\mathsf{type-A}$ problems, microscale modelling is required only in regions with physical defects or singularities (e.g., boundary layers, dislocations, shocks), while macroscale models are sufficient elsewhere. In contrast, $\mathsf{type-B}$ problems require microscale models throughout the domain, since a macroscale description may or may not be available, for instance, in complex system problems \cite{2007_roberts_general,2014_roberts_dynamical}, in heterogeneous problems \cite{2026_KUMARKARMAKAR_GPD}. The heterogeneous multiscale method ($\mathsf{HMM}$) and equation-free approaches fall into the $\mathsf{type-B}$ category, and therefore our proposed scheme in the computational domain is also of $\mathsf{type-B}$. The numerical examples in Subsection 4.1 exhibit boundary layers near the right and top boundaries in the physical domain; hence, a finescale is employed locally, while a coarse scale is used elsewhere, corresponding to $\mathsf{type-A}$. Consequently, the problems considered in this article are addressed using our proposed multiscale technique, which is applicable to both $\mathsf{type-A}$ and $\mathsf{type-B}$ settings.

\subsection{Non-axisymmetric diffusion in an annulus}

Various physical phenomena, such as mass transfer in chemical engineering and biology, as well as heat transfer in nuclear reactors, are effectively modelled using two-dimensional diffusion equations. In biological systems, multiscale techniques have been developed to simulate mass transport phenomena in tissues \cite{2021_yadav_heterogeneous,2017_islam_multiscale,2010_su_multi}. At the microscale, transport is governed by diffusion processes, while at the macroscale (tissue level), effective properties are inferred indirectly. Our method is capable of bridging these scales. Annular structures are commonly found in the membranes of living organisms. 
In cellular biology, annular regions are used to describe the diffusion of molecules or ions through the lipid bilayer of cell organelles or membranes \cite{1972_haydon_ion,2003_bockmann_effect}. Several critical physiological processes, such as nutrient transport and waste removal, depend on this mechanism \cite{1981_pardridge_transport,2007_das_multiscale}. Motivated by these real-world multiscale phenomena, we consider a non-axisymmetric diffusion problem in an annulus domain and solve it using the proposed patch dynamics scheme.

$\textbf{Physical Problem:}$
In this sub-section, we consider our physical domain $\Omega_p$ to be an annulus region bounded by the inner and outer circles $x^2+y^2$=1 and $x^2+y^2$=4, respectively. The problem in the physical domain $\Omega_p$ is equation \eqref{eqn:Physical_Problem}, where
\begin{equation}\label{eqn:Ann_Diff_Physical}
	\begin{split}
		D(x,y)=1, \hspace{0.2cm} v=\begin{bmatrix}
			0\hspace{0.3cm}0
		\end{bmatrix}^\top,\hspace{0.2cm} f=0, \hspace{0.2cm} g(x,y,t)=0.
	\end{split}
\end{equation}
The initial condition is
\begin{equation}
	u(x,y,0)=(\sqrt{x^2+y^2}-1)(2-\sqrt{x^2+y^2})\sin\left(\arctan\left(\frac{y}{x}\right)\right)\hspace{0.25cm} \text{in}\hspace{0.25cm} \Omega_p,
\end{equation} 
and the boundary conditions are homogeneous on both the inner and outer circles.

To map the physical annular domain onto a computational rectangular domain, an orthogonal transformation between the polar coordinates ($\xi$,$\eta$) and Cartesian coordinates ($x$,$y$) is employed:
\begin{equation}\label{eqn:Transformation_Annulus}
	x=\eta\cos\xi, \hspace{0.25cm} y=\eta\sin\xi, \hspace{0.25cm}\text{or} \hspace{0.25cm} \eta=\sqrt{x^2+y^2}, \hspace{0.25cm}\xi=\arctan \left(\frac{y}{x}\right).
\end{equation}

\begin{figure}%[h!]
	\centering
	
	\begin{tikzpicture}[scale=.38] 
		\draw[-latex] (-8,0)--(8,0) node[right] () {$x$};
		\draw[-latex] (0,-8)--(0,8) node[right] () {$y$};
		\polargrid[rmin=2,rmax=7,amin=0,amax=360]
		%----------------------------------
		\draw[fill=gray!50,opacity=2] (2.6704,0.842) -- (3.0519,0.9623) arc (17.5:27.5:3.2cm) -- (2.4836,1.2929) arc (27.5:17.5:2.8cm) -- cycle;
		
		\draw[fill=gray!50,opacity=2] (2.1449,1.7998) -- (2.4513,2.0569) arc (40:50:3.2cm) -- (1.7998,2.1449) arc (50:40:2.8cm) -- cycle;
		
		\draw[fill=gray!50,opacity=2] (1.2929
		,2.4836) -- (1.4776,2.8384) arc (62.5:72.5:3.2cm) -- (0.842,2.6704) arc (72.5:62.5:2.8cm) -- cycle;
		
		\draw[fill=gray!50,opacity=2] (0.244,2.7893) -- (0.2789,3.1878) arc (85:95:3.2cm) -- (-0.244,2.7893) arc (95:85:2.8cm) -- cycle;
		
		\draw[fill=gray!50,opacity=2] (-0.842,2.6704) -- (-0.9623,3.0519) arc (107.5:117.5:3.2cm) -- (-1.2929,2.4836) arc(117.5:107.5:2.8cm) -- cycle;
		
		\draw[fill=gray!50,opacity=2] (-1.7998,2.1449) -- (-2.0569,2.4513
		) arc (130:140:3.2cm) -- (-2.1449,1.7998) arc (140:130:2.8cm) -- cycle;
		
		\draw[fill=gray!50,opacity=2] (-2.4836,1.2929) -- (-2.8384,1.4776) arc (152.5:162.5:3.2cm) -- (-2.6704,0.842) arc (162.5:152.5:2.8cm) -- cycle;
		
		\draw[fill=gray!50,opacity=2] (-2.7893,0.244) -- (-3.1878,0.2789)  arc (175:185:3.2cm) -- (-2.7893,-0.244) arc (185:175:2.8cm) -- cycle;
		
		\draw[fill=gray!50,opacity=2] (-2.6704,-0.842) -- (-3.0519,-0.9623) arc (197.5:207.5:3.2cm) -- (-2.4836,-1.2929) arc (207.5:197.5:2.8cm) -- cycle;
		
		\draw[fill=gray!50,opacity=2] (-2.1449,-1.7998) -- (-2.4513,-2.0569) arc (220:230:3.2cm) -- (-1.7998,-2.1449) arc (230:220:2.8cm) -- cycle;
		
		\draw[fill=gray!50,opacity=2] (-1.2929,-2.4836) -- (-1.4776,-2.8384) arc (242.5:252.5:3.2cm) -- (-0.842,-2.6704) arc (252.5:242.5:2.8cm) -- cycle;
		
		\draw[fill=gray!50,opacity=2] (-0.244,-2.7893) -- (-0.2789,-3.1878) arc (265:275:3.2cm) -- (0.244,-2.7893) arc (275:265:2.8cm) -- cycle;
		
		\draw[fill=gray!50,opacity=2] (0.842,-2.6704) -- (0.9623,-3.0519) arc (287.5:297.5:3.2cm) -- (1.2929,-2.4836) arc (297.5:287.5:2.8cm) -- cycle;
		
		\draw[fill=gray!50,opacity=2] (1.7998,-2.1449) -- (2.0569,-2.4513) arc (310:320:3.2cm) -- (2.1449,-1.7998) arc (320:310:2.8cm) -- cycle;
		
		\draw[fill=gray!50,opacity=2] (2.4836,-1.2929) -- (2.8384,-1.4776) arc (332.5:342.5:3.2cm) -- (2.6704,-0.842) arc (342.5:332.5:2.8cm) -- cycle;
		
		%--------------------------------------------
		
		\draw[fill=blue!50,opacity=2] (3.6241,1.1427) -- (4.0056,1.263) arc (17.5:27.5:4.2cm) -- (3.3706,1.7546) arc (27.5:17.5:3.8cm) -- cycle;
		
		\draw[fill=blue!50,opacity=2] (2.911,2.4426) -- (3.2174,2.6997) arc (40:50:4.2cm) -- (2.4426,2.911) arc (50:40:3.8cm) -- cycle;
		
		\draw[fill=blue!50,opacity=2] (1.7546,3.3706) -- (1.9393,3.7254) arc (62.5:72.5:4.2cm) -- (1.1427,3.6241) arc (72.5:62.5:3.8cm) -- cycle;
		
		\draw[fill=blue!50,opacity=2] (0.3312,3.7855) -- (0.3661,4.184) arc (85:95:4.2cm) -- (-0.3312,3.7855) arc (95:85:3.8cm) -- cycle;
		
		\draw[fill=blue!50,opacity=2] (-1.1427,3.6241) -- (-1.263,4.0056) arc (107.5:117.5:4.2cm) -- (-1.7546,3.3706
		) arc(117.5:107.5:3.8cm) -- cycle;
		
		\draw[fill=blue!50,opacity=2] (-2.4426,2.911) -- (-2.6997,3.2174) arc (130:140:4.2cm) -- (-2.911,2.4426) arc (140:130:3.8cm) -- cycle;
		
		\draw[fill=blue!50,opacity=2] (-3.3706,1.7546) -- (-3.7254,1.9393) arc (152.5:162.5:4.2cm) -- (-3.6241,1.1427) arc (162.5:152.5:3.8cm) -- cycle;
		
		\draw[fill=blue!50,opacity=2] (-3.7855,0.3312) -- (-4.184,0.3661)  arc (175:185:4.2cm) -- (-3.7855,-0.3312
		) arc (185:175:3.8cm) -- cycle;
		
		\draw[fill=blue!50,opacity=2] (-3.6241,-1.1427) -- (-4.0056,-1.263) arc (197.5:207.5:4.2cm) -- (-3.3706
		,-1.7546) arc (207.5:197.5:3.8cm) -- cycle;
		
		\draw[fill=blue!50,opacity=2] (-2.911,-2.4426) -- (-3.2174,-2.6997) arc (220:230:4.2cm) -- (-2.4426,-2.911
		) arc (230:220:3.8cm) -- cycle;
		
		\draw[fill=blue!50,opacity=2] (-1.7546,-3.3706) -- (-1.9393,-3.7254) arc (242.5:252.5:4.2cm) -- (-1.1427
		,-3.6241) arc (252.5:242.5:3.8cm) -- cycle;
		
		\draw[fill=blue!50,opacity=2] (-0.3312,-3.7855) -- (-0.3661,-4.184) arc (265:275:4.2cm) -- (0.3312,-3.7855
		) arc (275:265:3.8cm) -- cycle;
		
		\draw[fill=blue!50,opacity=2] (1.1427,-3.6241) -- (1.263,-4.0056) arc (287.5:297.5:4.2cm) -- (1.7546,-3.3706) arc (297.5:287.5:3.8cm) -- cycle;
		
		\draw[fill=blue!50,opacity=2] (2.4426,-2.911) -- (2.6997,-3.2174) arc (310:320:4.2cm) -- (2.911,-2.4426) arc (320:310:3.8cm) -- cycle;
		
		\draw[fill=blue!50,opacity=2] (3.3706,-1.7546) -- (3.7254,-1.9393) arc (332.5:342.5:4.2cm) -- (3.6241,-1.1427) arc (342.5:332.5:3.8cm) -- cycle;
		%--------------------------------------------
		
		\draw[fill=red!50,opacity=2] (4.578,1.4434) -- (4.9593,1.5637) arc (17.5:27.5:5.2cm) -- (4.2577,2.2164) arc (27.5:17.5:4.8cm) -- cycle;
		
		\draw[fill=red!50,opacity=2] (3.677,3.0854) -- (3.9834,3.3425) arc (40:50:5.2cm) -- (3.0854,3.677) arc (50:40:4.8cm) -- cycle;
		
		\draw[fill=red!50,opacity=2] (2.2164,4.2577) -- (2.4011,4.6125) arc (62.5:72.5:5.2cm) -- (1.4434,4.5778) arc (72.5:62.5:4.8cm) -- cycle;
		
		\draw[fill=red!50,opacity=2] (0.4183,4.7817) -- (0.4532,5.1802) arc (85:95:5.2cm) -- (-0.4183,4.7817) arc (95:85:4.8cm) -- cycle;
		
		\draw[fill=red!50,opacity=2] (-1.4434,4.5778) -- (-1.5637,4.9593) arc (107.5:117.5:5.2cm) -- (-2.2164,4.2577) arc(117.5:107.5:4.8cm) -- cycle;
		
		\draw[fill=red!50,opacity=2] (-3.0854,3.677) -- (-3.3425,3.9834) arc (130:140:5.2cm) -- (-3.677,3.0854) arc (140:130:4.8cm) -- cycle;
		
		\draw[fill=red!50,opacity=2] (-4.2577,2.2164) -- (-4.6125,2.4011) arc (152.5:162.5:5.2cm) -- (-4.5778,1.4434) arc (162.5:152.5:4.8cm) -- cycle;
		
		\draw[fill=red!50,opacity=2] (-4.7817,0.4183) -- (-5.1802,0.4532)  arc (175:185:5.2cm) -- (-4.7817,-0.4183) arc (185:175:4.8cm) -- cycle;
		
		\draw[fill=red!50,opacity=2] (-4.58,-1.4434) -- (-4.9593,-1.5637) arc (197.5:207.5:5.2cm) -- (-4.2577,-2.2164) arc (207.5:197.5:4.8cm) -- cycle;
		
		\draw[fill=red!50,opacity=2] (-3.677,-3.0854) -- (-3.9834,-3.3425) arc (220:230:5.2cm) -- (-3.0854,-3.677) arc (230:220:4.8cm) -- cycle;
		
		\draw[fill=red!50,opacity=2] (-2.2164,-4.2577) -- (-2.4011,-4.6125) arc (242.5:252.5:5.2cm) -- (-1.4434,-4.5778) arc (252.5:242.5:4.8cm) -- cycle;
		
		\draw[fill=red!50,opacity=2] (-0.4183,-4.7817) -- (-0.4532,-5.1802) arc (265:275:5.2cm) -- (0.4183,-4.7817) arc (275:265:4.8cm) -- cycle;
		
		\draw[fill=red!50,opacity=2] (1.445,-4.5778) -- (1.5637,-4.9593) arc (287.5:297.5:5.2cm) -- (2.2164,-4.2577) arc (297.5:288.5:4.8cm) -- cycle;
		
		\draw[fill=red!50,opacity=2] (3.0854,-3.677) -- (3.3425,-3.9834) arc (310:320:5.2cm) -- (3.677,-3.0854) arc (320:310:4.8cm) -- cycle;
		
		\draw[fill=red!50,opacity=2] (4.2577,-2.2164) -- (4.6125,-2.4011) arc (332.5:342.5:5.2cm) --(4.5778,-1.4434) arc (342.5:332.5:4.8cm) -- cycle;
		
		%--------------------------------------------
		
		\draw[fill=brown!50,opacity=2] (5.5316,1.7441) -- (5.913,1.8644) arc (17.5:27.5:6.2cm) -- (5.1447,2.6781) arc (27.5:17.5:5.8cm) -- cycle;
		
		\draw[fill=brown!50,opacity=2] (4.4431,3.7282) -- (4.7495,3.9853) arc (40:50:6.2cm) -- (3.7282,4.4431) arc (50:40:5.8cm) -- cycle;
		
		\draw[fill=brown!50,opacity=2] (2.6781,5.1447) -- (2.8628,5.4995) arc (62.5:72.5:6.2cm) -- (1.7441,5.5316) arc (72.5:62.5:5.8cm) -- cycle;
		
		\draw[fill=brown!50,opacity=2] (0.5055,5.7779) -- (0.5404,6.1764) arc (85:95:6.2cm) -- (-0.5055,5.7779) arc (95:85:5.8cm) -- cycle;
		
		\draw[fill=brown!50,opacity=2] (-1.7441,5.5316) -- (-1.8644,5.913) arc (107.5:117.5:6.2cm) -- (-2.6781,5.1447) arc(117.5:107.5:5.8cm) -- cycle;
		
		\draw[fill=brown!50,opacity=2] (-3.7282,4.4431) -- (-3.9853,4.7495) arc (130:140:6.2cm) -- (-4.4431,3.7282) arc (140:130:5.8cm) -- cycle;
		
		\draw[fill=brown!50,opacity=2] (-5.1447,2.6781) -- (-5.4995,2.8628) arc (152.5:162.5:6.2cm) -- (-5.5316,1.7441) arc (162.5:152.5:5.8cm) -- cycle;
		
		\draw[fill=brown!50,opacity=2] (-5.7779,0.5055) -- (-6.1764,0.5404)  arc (175:185:6.2cm) -- (-5.7779,-0.5055) arc (185:175:5.8cm) -- cycle;
		
		\draw[fill=brown!50,opacity=2] (-5.5316,-1.7441) -- (-5.913,-1.8644) arc (197.5:207.5:6.2cm) -- (-5.1447,-2.6781) arc (207.5:197.5:5.8cm) -- cycle;
		
		\draw[fill=brown!50,opacity=2] (-4.4431,-3.7282) -- (-4.7495,-3.9853) arc (220:230:6.2cm) -- (-3.7282,-4.4431) arc (230:220:5.8cm) -- cycle;
		
		\draw[fill=brown!50,opacity=2] (-2.6781,-5.1447) -- (-2.8628,-5.4995) arc (242.5:252.5:6.2cm) -- (-1.7441,-5.5316) arc (252.5:242.5:5.8cm) -- cycle;
		
		\draw[fill=brown!50,opacity=2] (-0.5055,-5.7779) -- (-0.5404,-6.1764) arc (265:275:6.2cm) -- (0.5055,-5.7779) arc (275:265:5.8cm) -- cycle;
		
		\draw[fill=brown!50,opacity=2] (1.7441,-5.5316) -- (1.8644,-5.913) arc (287.5:297.5:6.2cm) -- (2.6781,-5.1447) arc (297.5:287.5:5.8cm) -- cycle;
		
		\draw[fill=brown!50,opacity=2] (3.7282,-4.4431) -- (3.9853,-4.7495) arc (310:320:6.2cm) -- (4.4431,-3.7282) arc (320:310:5.8cm) -- cycle;
		
		\draw[fill=brown!50,opacity=2] (5.1447,-2.6781) -- (5.4995,-2.8628) arc (332.5:342.5:6.2cm) --(5.5316,-1.7441) arc (342.5:332.5:5.8cm) -- cycle;
		
		\draw[very thick] (-0.8,0.8)  node {$\partial\Omega_{in}$};
		\draw[very thick] (-5.5,5.5)  node {$\partial\Omega_{out}$};
		\draw[very thick] (4.5,0.6)  node {$\Omega_p$};
		
		\draw [dashed] (2,0.1) -- (7,0.1);
		\draw [dashed] (2,-0.1) -- (7,-0.1);
		
		\filldraw[red] (2,0.1) circle (2pt) node[anchor=south]{$e$};
		\filldraw[black] (2,-0.1) circle (2pt) node[anchor=north]{$f$};
		\filldraw[brown] (7,0.1) circle (2pt) node[anchor=south]{$g$};
		\filldraw[blue] (7,-0.1) circle (2pt) node[anchor=north]{$h$};
		
		%			\draw[color=red] (1.4,0.5) circle [radius=0.4];
		%	\draw[color=red] (1.4,-0.5) circle [radius=0.4];
		%	\draw[color=red] (7,0.1) circle [radius=0.4];
		%	\draw[color=red] (7,-0.1) circle [radius=0.4];
		%--------------------------------------------	
		
		\draw[-latex] (10,-3.5)--(10,6) node[right] () {$\eta$};
		\draw[-latex] (10,-3.5)--(25.5,-3.5) node[right] () {$\xi$};
		\draw[thick] (10,-3.5) rectangle (24,3.5);
		\draw[thick] (10,-2.1) -- (24,-2.1);
		\draw[thick] (10,-0.7) -- (24,-0.7);
		\draw[thick] (10,0.7) -- (24,0.7);
		\draw[thick] (10,2.1) -- (24,2.1);

		\draw[thick] (10.875,-3.5) -- (10.875,3.5);
		\draw[thick] (11.75,-3.5) -- (11.75,3.5);
		\draw[thick] (12.625,-3.5) -- (12.625,3.5);
		\draw[thick] (13.5,-3.5) -- (13.5,3.5);
		\draw[thick] (14.375,-3.5) -- (14.375,3.5);
		\draw[thick] (15.25,-3.5) -- (15.25,3.5);
		\draw[thick] (16.125,-3.5) -- (16.125,3.5);
		\draw[thick] (17,-3.5) -- (17,3.5);
		\draw[thick] (17.875,-3.5) -- (17.875,3.5);
		\draw[thick] (18.75,-3.5) -- (18.75,3.5);
		\draw[thick] (19.625,-3.5) -- (19.625,3.5);
		\draw[thick] (20.5,-3.5) -- (20.5,3.5);
		\draw[thick] (21.375,-3.5) -- (21.375,3.5);
		\draw[thick] (22.25,-3.5) -- (22.25,3.5);
		\draw[thick] (23.125,-3.5) -- (23.125,3.5);
		
		\filldraw[fill=gray!50, draw=black] (10.625,-2.4) rectangle (11.125,-1.8);
		\filldraw[fill=gray!50, draw=black] (11.5,-2.4) rectangle (12,-1.8);
		\filldraw[fill=gray!50, draw=black] (12.375,-2.4) rectangle (12.875,-1.8);
		\filldraw[fill=gray!50, draw=black] (13.25,-2.4) rectangle (13.75,-1.8);
		\filldraw[fill=gray!50, draw=black] (14.125,-2.4) rectangle (14.625,-1.8);
		\filldraw[fill=gray!50, draw=black] (15,-2.4) rectangle (15.5,-1.8);
		\filldraw[fill=gray!50, draw=black] (15.875,-2.4) rectangle (16.375,-1.8);
		\filldraw[fill=gray!50, draw=black] (16.75,-2.4) rectangle (17.25,-1.8);
		\filldraw[fill=gray!50, draw=black] (17.625,-2.4) rectangle (18.125,-1.8);
		\filldraw[fill=gray!50, draw=black] (18.5,-2.4) rectangle (19,-1.8);
		\filldraw[fill=gray!50, draw=black] (19.375,-2.4) rectangle (19.875,-1.8);
		\filldraw[fill=gray!50, draw=black] (20.25,-2.4) rectangle (20.75,-1.8);
		\filldraw[fill=gray!50, draw=black] (21.125,-2.4) rectangle (21.625,-1.8);
		\filldraw[fill=gray!50, draw=black] (22,-2.4) rectangle (22.5,-1.8);
		\filldraw[fill=gray!50, draw=black] (22.875,-2.4) rectangle (23.375,-1.8);
		
		\filldraw[fill=blue!50, draw=black] (10.625,-1) rectangle (11.125,-0.4);
		\filldraw[fill=blue!50, draw=black] (11.5,-1) rectangle (12,-0.4);
		\filldraw[fill=blue!50, draw=black] (12.375,-1) rectangle (12.875,-0.4);
		\filldraw[fill=blue!50, draw=black] (13.25,-1) rectangle (13.75,-0.4);
		\filldraw[fill=blue!50, draw=black] (14.125,-1) rectangle (14.625,-0.4);
		\filldraw[fill=blue!50, draw=black] (15,-1) rectangle (15.5,-0.4);
		\filldraw[fill=blue!50, draw=black] (15.875,-1) rectangle (16.375,-0.4);
		\filldraw[fill=blue!50, draw=black] (16.75,-1) rectangle (17.25,-0.4);
		\filldraw[fill=blue!50, draw=black] (17.625,-1) rectangle (18.125,-0.4);
		\filldraw[fill=blue!50, draw=black] (18.5,-1) rectangle (19,-0.4);
		\filldraw[fill=blue!50, draw=black] (19.375,-1) rectangle (19.875,-0.4);
		\filldraw[fill=blue!50, draw=black] (20.25,-1) rectangle (20.75,-0.4);
		\filldraw[fill=blue!50, draw=black] (21.125,-1) rectangle (21.625,-0.4);
		\filldraw[fill=blue!50, draw=black] (22,-1) rectangle (22.5,-0.4);
		\filldraw[fill=blue!50, draw=black] (22.875,-1) rectangle (23.375,-0.4);
		
		\filldraw[fill=red!50, draw=black] (10.625,0.4) rectangle (11.125,1);
		\filldraw[fill=red!50, draw=black] (11.5,0.4) rectangle (12,1);
		\filldraw[fill=red!50, draw=black] (12.375,0.4) rectangle (12.875,1);
		\filldraw[fill=red!50, draw=black] (13.25,0.4) rectangle (13.75,1);
		\filldraw[fill=red!50, draw=black] (14.125,0.4) rectangle (14.625,1);
		\filldraw[fill=red!50, draw=black] (15,0.4) rectangle (15.5,1);
		\filldraw[fill=red!50, draw=black] (15.875,0.4) rectangle (16.375,1);
		\filldraw[fill=red!50, draw=black] (16.75,0.4) rectangle (17.25,1);
		\filldraw[fill=red!50, draw=black] (17.625,0.4) rectangle (18.125,1);
		\filldraw[fill=red!50, draw=black] (18.5,0.4) rectangle (19,1);
		\filldraw[fill=red!50, draw=black] (19.375,0.4) rectangle (19.875,1);
		\filldraw[fill=red!50, draw=black] (20.25,0.4) rectangle (20.75,1);
		\filldraw[fill=red!50, draw=black] (21.125,0.4) rectangle (21.625,1);
		\filldraw[fill=red!50, draw=black] (22,0.4) rectangle (22.5,1);
		\filldraw[fill=red!50, draw=black] (22.875,0.4) rectangle (23.375,1);
		
		\filldraw[fill=brown!50, draw=black] (10.625,1.8) rectangle (11.125,2.4);
		\filldraw[fill=brown!50, draw=black] (11.5,1.8) rectangle (12,2.4);
		\filldraw[fill=brown!50, draw=black] (12.375,1.8) rectangle (12.875,2.4);
		\filldraw[fill=brown!50, draw=black] (13.25,1.8) rectangle (13.75,2.4);
		\filldraw[fill=brown!50, draw=black] (14.125,1.8) rectangle (14.625,2.4);
		\filldraw[fill=brown!50, draw=black] (15,1.8) rectangle (15.5,2.4);
		\filldraw[fill=brown!50, draw=black] (15.875,1.8) rectangle (16.375,2.4);
		\filldraw[fill=brown!50, draw=black] (16.75,1.8) rectangle (17.25,2.4);
		\filldraw[fill=brown!50, draw=black] (17.625,1.8) rectangle (18.125,2.4);
		\filldraw[fill=brown!50, draw=black] (18.5,1.8) rectangle (19,2.4);
		\filldraw[fill=brown!50, draw=black] (19.375,1.8) rectangle (19.875,2.4);
		\filldraw[fill=brown!50, draw=black] (20.25,1.8) rectangle (20.75,2.4);
		\filldraw[fill=brown!50, draw=black] (21.125,1.8) rectangle (21.625,2.4);
		\filldraw[fill=brown!50, draw=black] (22,1.8) rectangle (22.5,2.4);
		\filldraw[fill=brown!50, draw=black] (22.875,1.8) rectangle (23.375,2.4);
		
		\filldraw[color=blue!60, fill=green!60, very thick](17,-3.5) circle (0.2);
		%	\filldraw[color=blue!60, fill=green!60, very thick](10,-3.5) circle (0.1);
		%	\filldraw[color=blue!60, fill=green!60, very thick](24,-3.5) circle (0.2);
		%	\filldraw[color=blue!60, fill=green!60, very thick](10,3.5) circle (0.2);
		
		\draw[very thick] (10,-4.85)  node {$\xi=0$};
		\draw[very thick] (17,-4.85)  node {$\xi=\pi$};
		\draw[very thick] (24,-4.85)  node {$\xi=2\pi$};
		\draw[very thick] (8.5,-3.5)  node {$\eta=1$};
		\draw[very thick] (8.5,3.5)  node {$\eta=2$};

		\filldraw[red] (10,-3.5) circle (2.5pt) node[anchor=north]{$e$};
		\filldraw[black] (24,-3.5) circle (2.5pt) node[anchor=north]{$f$};
		\filldraw[blue] (24,3.5) circle (2.5pt) node[anchor=south]{$h$};
		\filldraw[brown] (10,3.5) circle (2.5pt) node[anchor=south]{$g$};

	\end{tikzpicture}  
	\caption{Left: Physical domain \& Right: Computational domain} 
	\label{fig:Annulus_Phy_Comp}
\end{figure}	

A schematic view of the physical annular domain and the corresponding computational rectangular domain is shown in Figure \ref{fig:Annulus_Phy_Comp}. In the physical domain $\Omega_p$, the patches are non-rectangular in shape. However, under the transformation defined in equation \eqref{eqn:Transformation_Annulus}, these non-rectangular patches correspond to rectangular patches in the computational domain $\Omega_c=[0,2\pi)\times[1,2]$. Here, $J=r>1$ throughout the entire domain. The transformation \eqref{eqn:Transformation_Annulus} is smooth, one-to-one, and $1<J<2$, which shows that the mapping is non-singular.

$\textbf{Computational Problem:}$

Under the transformation \eqref{eqn:Transformation_Annulus}, the physical problem \eqref{eqn:Physical_Problem} together with \eqref{eqn:Ann_Diff_Physical} is reduced to the computational problem \eqref{eqn:Computational_Problem}, where
\begin{equation}\label{eqn:Ann_Diff_Computational}
	\begin{split}
		&\alpha(\xi,\eta)=\frac{1}{\eta^2}, \hspace{0.2cm} \beta(\xi,\eta)=0,\hspace{0.2cm}
		\gamma(\xi,\eta)=1, \hspace{0.2cm}\\
		&\nu(\xi,\eta)=0\hspace{0.2cm}
		\omega(\xi,\eta)=-\frac{1}{\eta}\hspace{0.2cm}
		\phi(\xi,\eta)=0, \hspace{0.2cm} g(\xi,\eta,t)=0. 
	\end{split}
\end{equation}
The corresponding initial and boundary conditions in $\Omega_c$ are:
\begin{equation}\label{eqn:Ann_C_IC_BCs}
	\begin{split}
		&\mathsf{IC:}\hspace{0.25cm}u(\xi,\eta,0)=(\eta-1)(2-\eta)\sin\xi,\hspace{0.25cm} \text{in}\hspace{0.25cm} \Omega_c,\\
		&\mathsf{BCs:}\hspace{0.25cm}u(0,\eta,t)=u(2\pi,\eta,t), \hspace{0.25cm} \frac{\partial u}{\partial\xi}(0,\eta,t)=\frac{\partial u}{\partial\xi}(2\pi,\eta,t),\\
		&\hspace{1.05cm}  u(\xi,1,t)=0,\hspace{0.25cm}u(\xi,2,t)=0,\hspace{0.25cm} \text{on}\hspace{0.2cm} \partial\Omega_c \hspace{0.2cm} \text{and}\hspace{0.2cm} t\in[0,0.2].
	\end{split}
\end{equation}

The analytical solution of the computational problem is given by 
\begin{equation}\label{eqn:Annulus_AS}
	u_a(\xi,\eta,t)=\frac{\pi^2}{2}\sum_{m=1}^{\infty}\left[\int_{1}^{2}\psi(\psi-1)(2-\psi)Z_1(\mu_{1,m}\psi) d\psi\right]B_{1,m}Z_1(\mu_{1,m}\eta)\sin\xi \exp(-\mu_{1,m}^2t),
\end{equation}
where 
\begin{center}
	$B_{1,m}=\frac{\mu_{1,m}^2 J_1^2(2\mu_{1,m})}{J_1^2(\mu_{1,m})-J_1^2(2\mu_{1,m})}$,
\end{center}
$Z_1(\mu_{1,m}\eta)=J_1(\mu_{1,m})Y_1(\mu_{1,m}\eta)-Y_1(\mu_{1,m}) J_1(\mu_{1,m}\eta)$, $J_1$ and $Y_1$ are the Bessel's functions of first kind and second kind respectively. $\mu_{1,m}$ are the positive roots of the transcendental equation $J_1(\mu)Y_1(2\mu)-Y_1(\mu)J_1(2\mu)=0$. 

Using the method of lines for spatial discretisation, the central difference scheme is employed for both the radial and azimuthal derivatives in equation \eqref{eqn:Computational_Problem} together with \eqref{eqn:Ann_Diff_Computational}. Then we have 
\begin{equation}\label{eqn:MOL_Diff_Annu}
	\begin{split}
		\frac{d}{dt}u_{i,j}=&\frac{1}{\delta \eta^2}(u_{i,j+1}-2u_{i,j}+u_{i,j-1})+\frac{1}{2\eta_j\hspace{0.05cm}\delta \eta}(u_{i,j+1}-u_{i,j-1})\\
		&\frac{1}{\eta_j^2\hspace{0.05cm}\delta\xi^2}(u_{i+1,j}-2u_{i,j}+u_{i-1,j}),
	\end{split}
\end{equation}
where $\delta\xi$ and $\delta\eta$ are the spatial nano-steps along $\xi$- and $\eta$-directions, respectively. %$u_{i,j}(t)$ is a $(\operatorname{n\theta}+1)\times(\operatorname{nr}+1)$-dimensional vector.
$\operatorname{n}$ denotes the number of spatial nano grid steps along both $\xi$- and $\eta$-directions, respectively. Small values of $\delta\xi$ and $\delta\eta$ show a system of stiff equations \eqref{eqn:MOL_Diff_Annu} \cite{2003_gear_projective,2003_gear_telescopic,2007_lee_second}. %Classical numerical integration techniques are not appropriate to solve stiff $\mathsf{ODEs}$ in the entire time interval $[0, T]$ because of the time step restriction. So we need an efficient multiscale technique to deal with this kind of problem. 
We apply our proposed patch dynamics scheme on the stiff problem \eqref{eqn:MOL_Diff_Annu} to find a system-level solution in the large domain over a long time.

In order to solve the problem \eqref{eqn:MOL_Diff_Annu} using the proposed scheme, we consider $h_\xi=h_\eta=0.001$ and a time-stepper of size $\tau=1e{-6}$. At the microscopic level, each patch is discretised with $\operatorname{n}=20$ spatial nano grid steps in both azimuthal and radial directions, and the time-stepper is discretised using $\operatorname{n_t}=1500$ nano time steps. The values of $\delta\xi$ and $\delta\eta$ are same and equal to $5\mathrm{e}{-5}$, which shows that equation \eqref{eqn:MOL_Diff_Annu} is a system of stiff $\mathsf{ODEs}$. An explicit forward Euler scheme is used for the temporal derivative in the micro simulation. The trapezoidal composite rule is used in the restriction operator \eqref{eqn:Restriction} to restrict the microscopic values in the patch to the macroscopic value. After obtaining the patch dynamics solution in the computational domain $\Omega_c$, the inverse transformation \eqref{eqn:Transformation_Annulus} is applied to recover the patch dynamics solution in the physical domain $\Omega_p$.

\begin{table}%[h!]
	%	\setstretch{1.2}
%	\begin{center}
	\caption{%\setstretch{1.2}
		The second column presents the maximum percentage errors in the patch dynamics solution ($U$) relative to the analytic solution ($u_a$) for the problem \eqref{eqn:Physical_Problem} with \eqref{eqn:Ann_Diff_Physical}, corresponding to various micro-simulation times ($\tau$) listed in the first column.}
				\centering
			\begin{tabular}{cc}
				%				\hline
				%				\multicolumn{1}{|c|}{} & \multicolumn{2}{c|}{Maximum percentage errors for $(\operatorname{N_x}+1)\times(\operatorname{N_y}+1), \operatorname{N_t}$} \\
				\hline
				$\tau$ & Max. Per. Error\\ %$\left\|\frac{U-u_a}{u_a}\right\|_\infty\times100\%$ \\
				\hline
				$1\mathrm{e}{-6}$&0.55\\
				$2\mathrm{e}{-6}$&0.53\\
				$4\mathrm{e}{-6}$&0.52\\
				$6\mathrm{e}{-6}$&0.51\\
				$8\mathrm{e}{-6}$&0.51\\
				$10\mathrm{e}{-6}$&0.51\\
				\hline
		\end{tabular}
%	\end{center}
	\label{table:relaxation_time_lambda0_Diff}
\end{table}

For problem \eqref{eqn:Physical_Problem} with \eqref{eqn:Ann_Diff_Physical}, Table \ref{table:relaxation_time_lambda0_Diff} compares the accuracy of the patch dynamics solutions for different micro time duration $\tau$. For this experiment, we set $\operatorname{N_x}=16$, $\operatorname{N_y}=10$, $\operatorname{N_t}=500$ and used same microscale discretisation as described above. During this experiment, we varied only the micro time step $\tau$ while keeping all other parameters fixed. As shown in Table \ref{table:relaxation_time_lambda0_Diff}, increasing the micro time step ($\tau$) results in no noticeable change in the maximum percentage error. Based on the above discussions, we chose the micro time step $\tau=1e-6$ for all subsequent experiments on the non-axisymmetric diffusion in an annulus region in Sub-section 4.2.

%\begin{table}%[h!]
%	\setstretch{1.2}
%	\begin{center}
	%		\scalebox{0.7}{
		%			\begin{tabular}{|p{3cm}|p{2cm}|p{2cm}|p{2cm}|p{2cm}|}
			%				\hline
			%				\multicolumn{1}{|c|}{} & \multicolumn{4}{c|}{Maximum percentage errors at time $t$}\\
			%				\hline
			%				($\operatorname{N_x}, \operatorname{N_y}, \operatorname{N_t}$)&$t=0.25$&$t=0.5$&$t=0.75$&$t=1$\\
			%				\hline
			%				(16, 10, 500)&0.6824&1.3115&1.9362&2.5569\\
			%				(24, 15, 1200)&0.2582&0.4961&0.7333&0.9700\\
			%				(32, 20, 2000)&0.2064&0.4042&0.6014&0.7983\\
			%				(40, 25, 3200)&0.1387&0.2745&0.4100&0.5454\\
			%				(48, 30, 4500)&0.1169&0.2341&0.3510&0.4679\\
			%				(56, 35, 6200)&0.0952&0.1927&0.2900&0.3872\\
			%				(64, 40, 8000)&0.0854&0.1743&0.2631&0.3518\\
			%				\hline
			%		\end{tabular}}
	%	\end{center}
%\vspace{0.2cm}
%	\caption{\setstretch{1.2}The percentage errors in the patch dynamics solution of the problem \eqref{eqn:Ann_Diff_Physical} are presented for various spatial resolutions and at different time instances.}
%	\label{table:P2_PerErr}
%\end{table}

\begin{table}%[h!]
%	\begin{center}
	\caption{The maximum relative errors in the patch dynamics solution of the problem \eqref{eqn:Physical_Problem} together with \eqref{eqn:Ann_Diff_Physical} are presented for various spatial resolutions and at different time instances.}
	\centering
		\scalebox{1}{
			\begin{tabular}{p{2.5cm}p{1.5cm}p{1.5cm}p{1.5cm}p{1.5cm}}
				%				\hline
				%				\multicolumn{1}{|c|}{} & \multicolumn{4}{c|}{Maximum percentage errors at time $t$}\\
				\hline
				($\operatorname{N_x}, \operatorname{N_y}, \operatorname{N_t}$)&$t=0.05$&$t=0.1$&$t=0.15$&$t=0.2$\\
				\hline
				(16, 10, 500)&1.40$\mathrm{e}{-3}$&2.85$\mathrm{e}{-3}$&4.24$\mathrm{e}{-3}$&5.55$\mathrm{e}{-3}$\\
				(24, 15, 1200)&5.40$\mathrm{e}{-4}$&1.09$\mathrm{e}{-3}$&1.61$\mathrm{e}{-3}$&2.10$\mathrm{e}{-3}$\\
				(32, 20, 2000)&3.75$\mathrm{e}{-4}$&8.23$\mathrm{e}{-4}$&1.26$\mathrm{e}{-3}$&1.66$\mathrm{e}{-3}$\\
				(40, 25, 3200)&2.31$\mathrm{e}{-4}$&5.38$\mathrm{e}{-4}$&8.33$\mathrm{e}{-4}$&1.11$\mathrm{e}{-3}$\\
				(48, 30, 4500)&1.74$\mathrm{e}{-4}$&4.38$\mathrm{e}{-4}$&6.91$\mathrm{e}{-4}$&9.32$\mathrm{e}{-4}$\\
				(56, 35, 6200)&1.41$\mathrm{e}{-4}$&3.46$\mathrm{e}{-4}$&5.56$\mathrm{e}{-4}$&7.56$\mathrm{e}{-4}$\\
				(64, 40, 8000)&1.92$\mathrm{e}{-4}$&3.02$\mathrm{e}{-4}$&4.93$\mathrm{e}{-4}$&6.75$\mathrm{e}{-4}$\\
				\hline
		\end{tabular}}
%	\end{center}
	\label{table:P2_PerErr}
\end{table}

Table \ref{table:P2_PerErr} presents the maximum relative errors in the patch dynamics solution with respect to the analytical solution \eqref{eqn:Annulus_AS} of the physical problem \eqref{eqn:Physical_Problem} together with \eqref{eqn:Ann_Diff_Physical}, evaluated at various grid resolutions and time instances $t=0.05$, 0.1, 0.15 and 0.2. The results demonstrate that grid refinement leads to improved accuracy, indicating the convergence of the patch dynamics solution to the analytical solution. Since the physical problem \eqref{eqn:Physical_Problem} together with \eqref{eqn:Ann_Diff_Physical} is a diffusion equation with homogeneous boundary conditions, the solution \eqref{eqn:Annulus_AS} naturally decays over time. Consequently, the relative error tends to increase as the solution magnitude decreases. Overall, the results in Table \ref{table:P2_PerErr} show that the patch dynamics solutions are in good agreement with the analytical solution. 

\begin{figure}%[h!]
	\begin{subfigure}{.45\textwidth}
		\centering
		% include first image
		\includegraphics[width=1\linewidth]{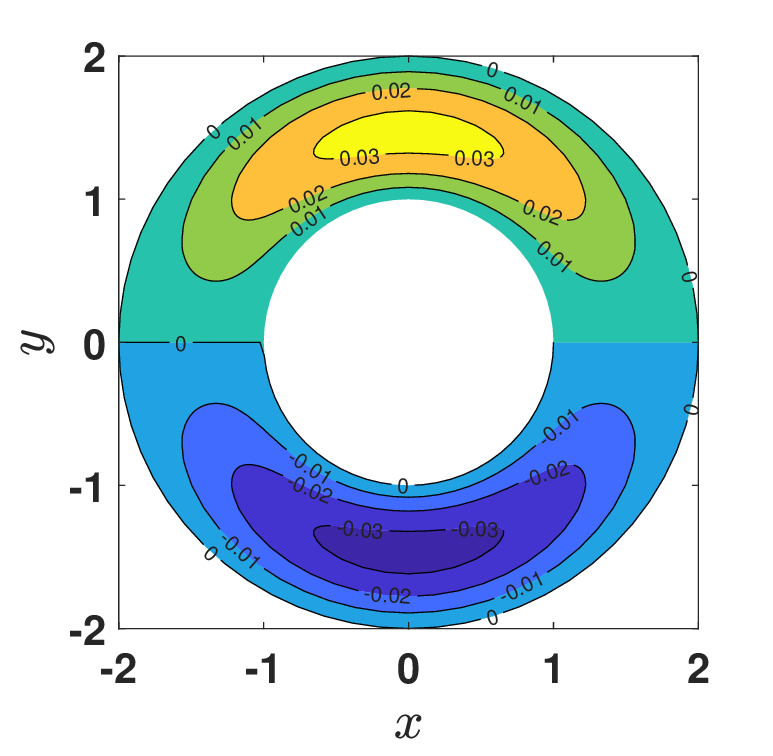}  
		\caption{Analytical solution}
		\label{fig:P2_Soln_Contour_AS}
	\end{subfigure}
	\begin{subfigure}{.45\textwidth}
		\centering
		% include second image
		\includegraphics[width=1\linewidth]{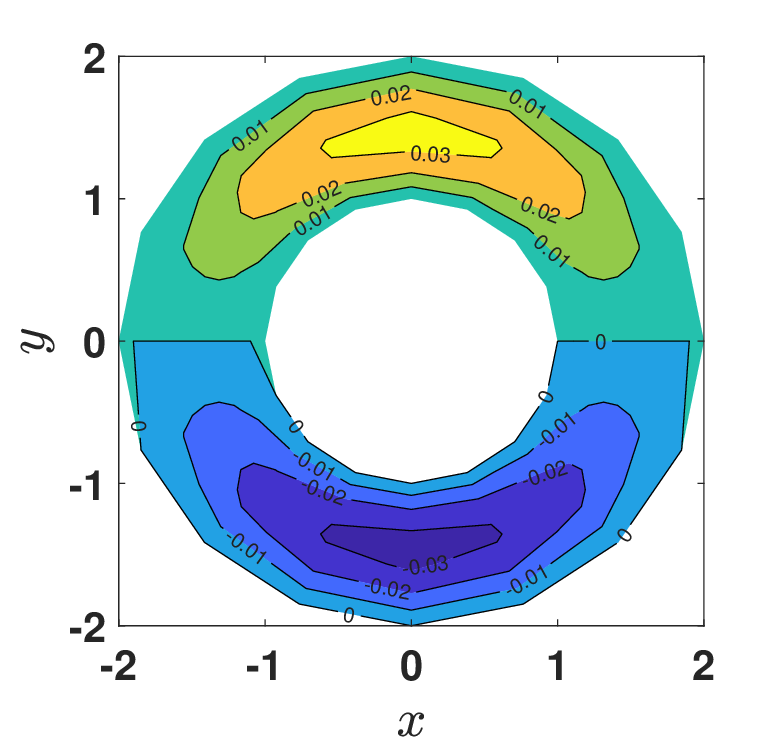}  
		\caption{16$\times$10 grid}
		\label{fig:P2_Soln_Contour_17_11}
	\end{subfigure}
	
	%		\newline
	
	\begin{subfigure}{.45\textwidth}
		\centering
		% include third image
		\includegraphics[width=1\linewidth]{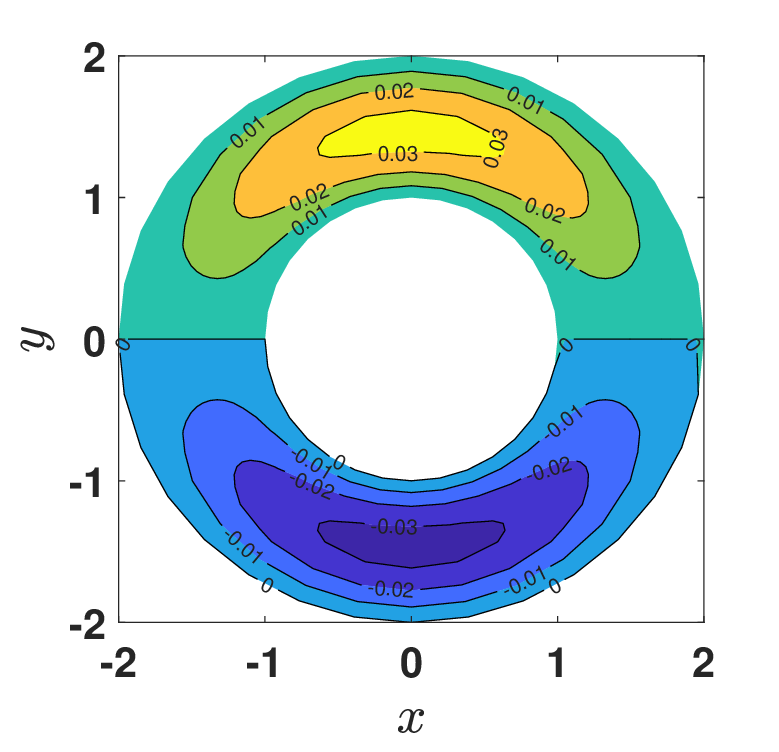}  
		\caption{32$\times$20 grid}
		\label{fig:P2_Soln_Contour_33_21}
	\end{subfigure}
	\begin{subfigure}{.45\textwidth}
		\centering
		% include fourth image
		\includegraphics[width=1\linewidth]{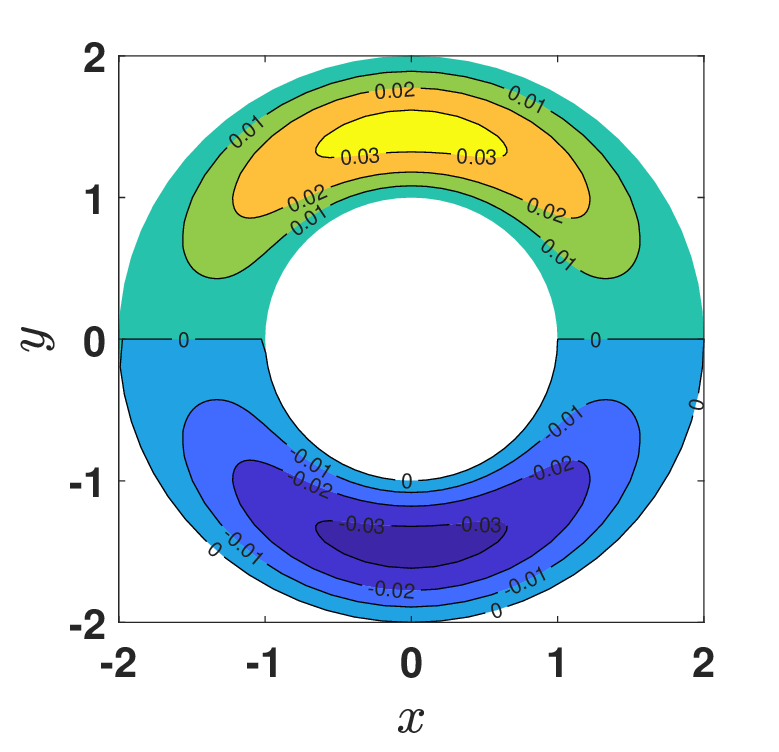}  
		\caption{64$\times$40 grid}
		\label{fig:P2_Soln_Contour_65_41}
	\end{subfigure}
	\caption{Analytical solution and patch dynamics solutions of the problem \eqref{eqn:Physical_Problem} together with \eqref{eqn:Ann_Diff_Physical} are presented for different grid resolutions in the physical domain $\Omega_p$ at final time $T=0.2$.}
	\label{fig:P2_Soln_Contour}
\end{figure}

Figure \ref{fig:P2_Soln_Contour} displays contour plots of the analytical as well as patch dynamics solutions for three different grid resolutions at final time $T=0.2$. All patch dynamics solutions exhibit good agreement with the analytical solution. However, in Figure \ref{fig:P2_Soln_Contour_17_11}, lower macro grid resolution causes the curved contour lines to appear less smooth.

\begin{figure}
	\begin{subfigure}{.5\textwidth}
		\centering
		% include first image
		\includegraphics[width=1\linewidth]{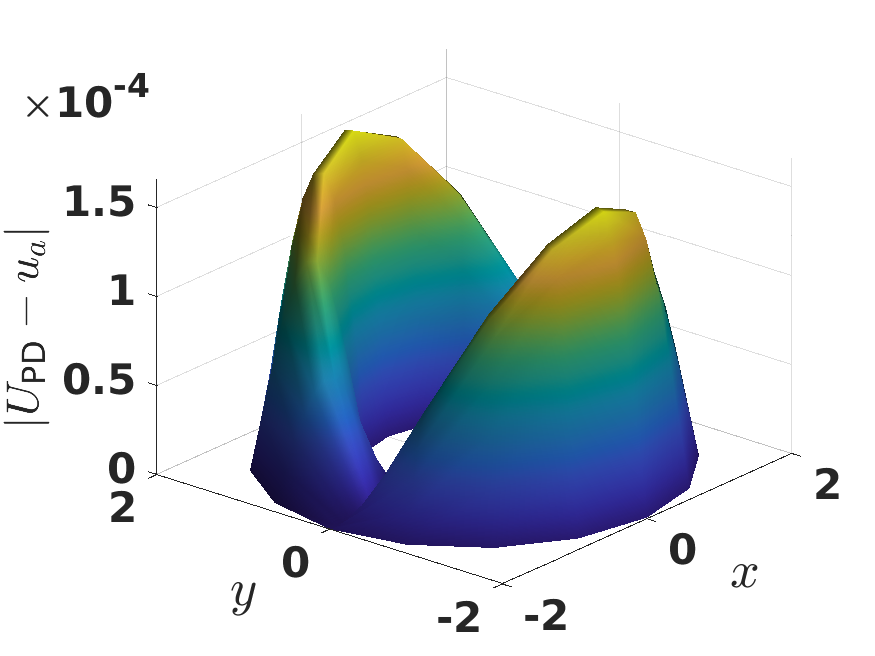}  
		\caption{16$\times$10 grid}
		\label{fig:P2_fig12}
	\end{subfigure}
	\begin{subfigure}{.5\textwidth}
		\centering
		% include second image
		\includegraphics[width=1\linewidth]{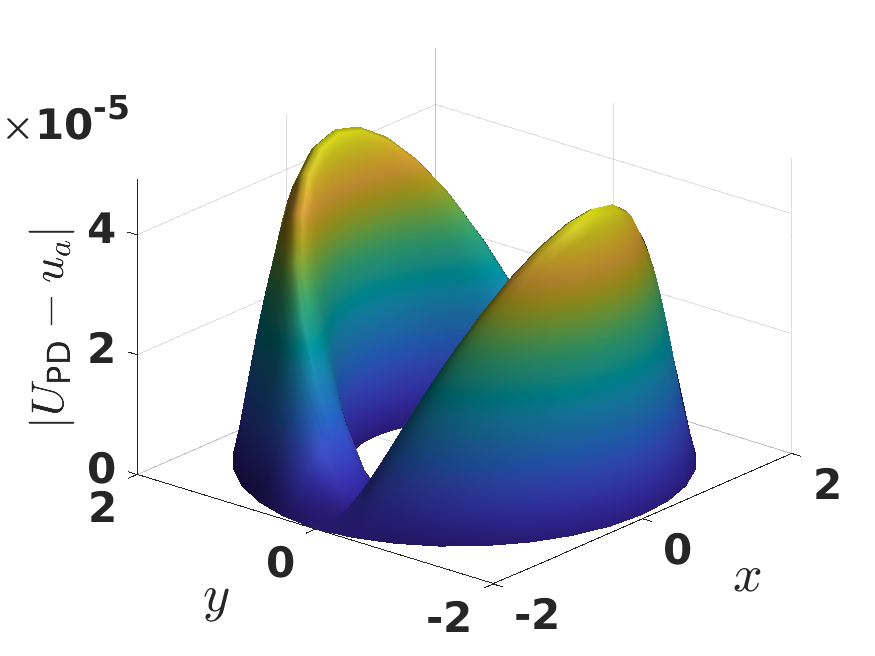}  
		\caption{32$\times$20 grid}
		\label{fig:P2_fig13}
	\end{subfigure}
	
	%		\newline
	\centering
	\begin{subfigure}{.5\textwidth}
		\centering
		% include third image
		\includegraphics[width=1\linewidth]{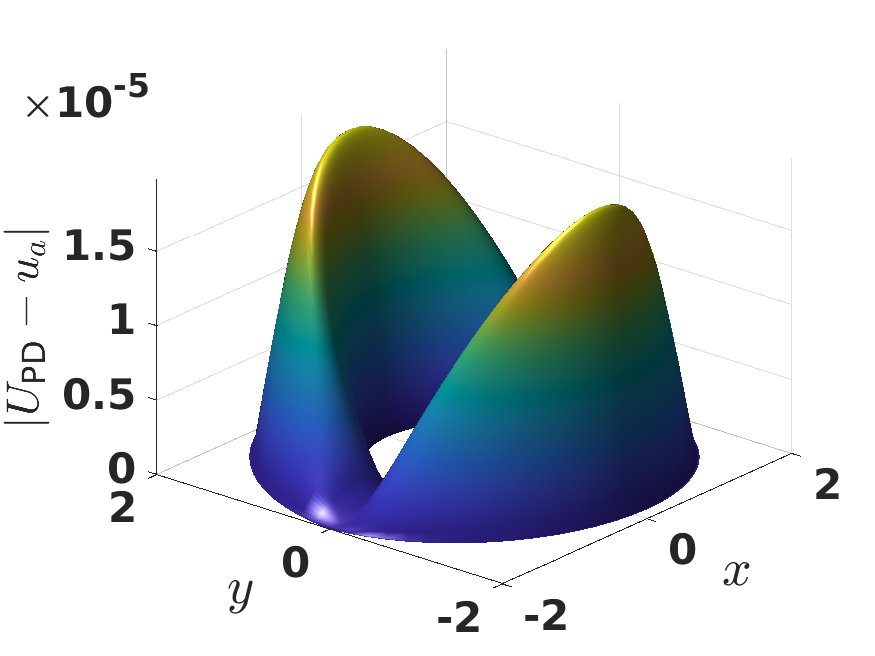}  
		\caption{64$\times$40 grid}
		\label{fig:}
	\end{subfigure}
	\caption{Absolute errors in the patch dynamics solution of the problem \eqref{eqn:Physical_Problem} together with \eqref{eqn:Ann_Diff_Physical} in the whole domain with various spatial resolutions are presented at the final time $T=0.2$. $U_\mathsf{PD}$ and $u_a$ denote patch dynamics and analytical solutions of the physical problem \eqref{eqn:Physical_Problem} together with \eqref{eqn:Ann_Diff_Physical}.}
	\label{fig:P2_AbsError_Surface}
\end{figure}

Figure \ref{fig:P2_AbsError_Surface} presents surface plots of the absolute errors in the patch dynamics solutions for three grid resolutions $16\times10$, $32\times20$, and $64\times40$ at final time $T=0.2$. The corresponding maximum percentage errors are 1.66$\mathrm{e}{-4}$, 4.93$\mathrm{e}{-5}$, and 1.98$\mathrm{e}{-5}$, respectively. The results for problem \eqref{eqn:Physical_Problem} demonstrate that the proposed patch dynamics scheme effectively handles periodic boundary conditions at the macroscopic level in the computational domain. Furthermore, these findings confirm that by using a body-fitted orthogonal curvilinear coordinate system in combination with an efficient patch dynamics scheme, one can accurately predict the system-level behaviour in a non-rectangular physical domain $\Omega_p$ by performing microscopic simulation in a rectangular computational domain $\Omega_c$.

%%%%%%%%%%%%%%%%%%%%%%%%%%%%%%%%%%%%%%%%%%%%%%%%%%%%%%%%%%%%%

\section{Conclusion}

So far, patch dynamics schemes have mainly been developed on rectangular domains using uniform rectangular patches in shape and size. However, considering the complexities arising from both physical and geometrical aspects of various real-life problems, it is evident that domains may not necessarily be rectangular. An explicit patch dynamics scheme is proposed in this study to address the challenges posed by two-dimensional, unsteady, linear, heterogeneous convection-diffusion-reaction ($\mathsf{CDR}$) equations, where non-uniform grids are used in order to determine the system-level behaviour in large non-rectangular domains over a long time. The shape and size of each patch are determined based on the physical and geometrical complexities of the problem, and generalised orthogonal curvilinear coordinates are employed to handle such complexities. Based on such complexities, two different types of numerical problems are chosen for validation. In the first problem, a convection-dominated two-dimensional heterogeneous $\mathsf{CDR}$ equation is considered over a rectangular domain with high gradients (or boundary layer) near the right and top boundaries.
%Non-homogeneous Dirichlet boundary conditions are considered, exhibiting exponential growth along the $x$- and $y$- directions. In the first and second parts of the problem, constant and variable diffusivities are considered, and their solutions are verified with the analytical solution and some existing numerical solutions. 
%The grid independence of both patch dynamics solutions is demonstrated on $21\times 21$ macro grid over $\operatorname{Nt}=8500$ macro time steps.
In the second one, a two-dimensional unsteady non-axisymmetric diffusion problem is solved in an annulus using the body-fitted orthogonal curvilinear coordinate system. There is an excellent agreement between the solutions obtained through the proposed scheme and the existing findings. Based on the outcomes of the proposed scheme, the following conclusions are drawn:

\begin{enumerate}
	\item A methodical representation of a patch dynamics scheme is proposed for general unsteady, linear, heterogeneous convection-diffusion-reaction ($\mathsf{CDR}$) equations on a generalised orthogonal curvilinear coordinate system in a two-dimensional domain.
	\item The patch dynamics scheme is tailored to solve multiscale problems on non-uniform grids and non-rectangular domains, accommodating non-uniform and non-rectangular patch configurations within the physical domain. 
	\item The proposed scheme is capable of handling convection-dominated unsteady heterogeneous $\mathsf{CDR}$ equations on a two-dimensional orthogonal curvilinear coordinate system, where the diffusion tensor and the convection velocity depend explicitly on the spatial variables. This scheme also efficiently handles space-time dependent source terms.
	\item A stretched grid is employed for such problems where high-gradient (or boundary layer) regions are observed. Stretched grids with a stretching ratio, $\lambda=0.1$, provide a better accuracy compared to unstretched grids (i.e., uniform grids where $\lambda=0$) and other stretching ratios ($\lambda\ne 0.1$). Stretched grids notably reduce errors substantially in the higher-gradient regions (or in boundary layer regions).
	\item Grid independence of the patch dynamics solutions for both constant and variable diffusivities in the heterogeneous $\mathsf{CDR}$ equations is performed, and it is observed that a $20\times 20$ macro grid with $\operatorname{N_t}=8500$ macro time steps provides grid independent results with sufficient accuracy. %In this resolution, a balance is achieved between computational efficiency and solution accuracy.
	\item Overall, considering accuracy, computational time and memory usage, the proposed patch dynamics scheme performs much better compared to the full-domain simulation. 
	%	\item The problems considered in this study are addressed using our proposed multiscale technique, which is applicable to both $\mathsf{type-A}$ and $\mathsf{type-B}$ settings.
	\item The proposed patch dynamics scheme can efficiently handle a non-axisymmetric diffusion problem in an annulus with periodic as well as Dirichlet boundary conditions.
	\item Numerical results show that the patch dynamics solution converges to the true solution as the macroscopic grids and time levels are refined.
\end{enumerate}

As an initial development of the patch dynamics scheme in two-dimensional curvilinear coordinates, this article focuses exclusively on linear, unsteady, heterogeneous $\mathsf{CDR}$ equations. Further research is required to extend the framework to heterogeneous problems with highly oscillatory coefficients.

%\begin{thebibliography}{00}
%
%%% For numbered reference style
%%% \bibitem{label}
%%% Text of bibliographic item
%
%%\bibitem{lamport94}
%%  Leslie Lamport,
%%  \textit{\LaTeX: a document preparation system},
%%  Addison Wesley, Massachusetts,
%%  2nd edition,
%%  1994.
%
%\end{thebibliography}

\bibliographystyle{elsarticle-num}
\bibliography{Work2}
\end{document}